\documentclass{amsart}

\usepackage{amssymb,amsfonts, latexsym,amsmath,hhline,array,longtable}
\usepackage{pdfsync,color, comment,colortbl, mathrsfs,stmaryrd,cite,graphicx,yfonts,amsthm}
\usepackage[colorlinks=true, citecolor=blue, linkcolor=blue, urlcolor=blue]{hyperref}

\input xy
\xyoption{all}

\def\grphp#1{$\xymatrix@R=10pt@C=10pt@M=0pt@L=2pt{#1}$}

\newcommand{\cal}{\mathcal}

\newtheorem{formula}{}[section]
\newtheorem{definition}[formula]{Definition}
\newtheorem{corollary}[formula]{Corollary}
\newtheorem{remark}[formula]{Remark}
\newtheorem{lemma}[formula]{Lemma}
\newtheorem{theorem}[formula]{Theorem}
\newcounter{ccntr}
\newtheorem{claim}[ccntr]{Claim}

\def\thrm{\begin{theorem}}
\def\thrml#1{\begin{theorem}\label{#1}}
\def\ethrm{\end{theorem}}
\def\rmrk{\begin{remark}}
\def\rmrkl#1{\begin{remark}\label{#1}}
\def\ermrk{\end{remark}}
\def\dfntn{\begin{definition}}
\def\dfntnl#1{\begin{definition}\label{#1}}
\def\edfntn{\end{definition}}
\def\nmrt{\begin{enumerate}}
\def\enmrt{\end{enumerate}}
\def\tm#1{\item[{\rm (#1)}]}
\def\qtnl#1{\begin{equation}\label{#1}}
\def\eqtn{\end{equation}}
\def\lmm{\begin{lemma}}
\def\lmml#1{\begin{lemma}\label{#1}}
\def\elmm{\end{lemma}}
\def\crllr{\begin{corollary}}
\def\crllrl#1{\begin{corollary}\label{#1}}
\def\ecrllr{\end{corollary}}
\def\css{\begin{cases}}
\def\ecss{\end{cases}}
\def\prf{\begin{proof}}
\def\eprf{\end{proof}}
\def\clm{\begin{claim}}
\def\clml#1{\begin{claim}\label{#1}}
\def\eclm{\end{claim}}

\def\cA{{\cal A}}

\def\cG{{\cal G}}

\def\cR{{\cal R}}

\def\cS{{\cal S}}

\def\cX{{\cal X}}

\def\mN{{\mathbb N}}

\def\mZ{{\mathbb Z}}

\def\fB{{\mathfrak B}}
\def\fC{{\mathfrak C}}
\def\fc{{\mathfrak c}}

\def\fn{{\mathfrak n}}
\def\fM{{\mathfrak M}}
\def\fm{{\mathfrak m}}

\def\fs{{\mathfrak s}}
\def\fS{{\mathfrak S}}
\def\fT{{\mathfrak T}}
\def\ft{{\mathfrak t}}

\def\fX{{\mathfrak X}}

\def\sC{{\mathscr{C}}}
\def\sR{{\mathscr{R}}}
\def\sX{{\mathscr{X}}}

\DeclareMathOperator{\aut}{Aut}

\DeclareMathOperator{\cay}{Cay}

\DeclareMathOperator{\cyc}{cyc}

\DeclareMathOperator{\id}{id}
\DeclareMathOperator{\im}{im}

\DeclareMathOperator{\inv}{Inv}

\DeclareMathOperator{\iso}{Iso}

\DeclareMathOperator{\mult}{\fM}
\DeclareMathOperator{\omult}{\fC}

\DeclareMathOperator{\pr}{pr}
\DeclareMathOperator{\PSL}{PSL}

\DeclareMathOperator{\rad}{rad}
\DeclareMathOperator{\res}{Res}
\DeclareMathOperator{\rk}{rk}

\DeclareMathOperator{\Span}{Span}

\DeclareMathOperator{\sym}{Sym}
\DeclareMathOperator{\WL}{WL}

\def\bone{{\bf 1}}

\def\grp#1{\langle {#1}\rangle}

\def\phmaa#1{{\phantom{x}\hspace{-2mm}^{#1}}}
\def\qaq{\quad\text{and}\quad}
\def\qoq{\quad\text{or}\quad}
\def\ov{\overline}
\def\und#1{{\underline{#1}}}
\def\wh{\widehat}

\begin{document}

\title[A lower bound for the $\WL$-dimension of circulants]{A lower bound for the Weisfeiler-Leman dimension of circulant graphs}
\author{Yulai Wu}
\address{School of Mathematics and Statistics, Hainan University, Haikou, China}
\email{wuyl@hainanu.edu.cn}
\thanks{Partly supported by Hainan Province Natural Science Foundation of China, grant No. 120RC452.}
\author{Qing Ren*}
\address{Guangxi Technological College of Machinery and Electricity, Nanning, China}
\thanks{Partly supported by National Science Foundation of China, grant No. 12371019}
\email{renqing@gxcme.edu.cn}
\author{Ilia Ponomarenko}
\address{Steklov Institute of Mathematics at St. Petersburg, Russia}
\email{inp@pdmi.ras.ru}
\thanks{}
\thanks{$^*$Corresponding author: Qing Ren}
\date{}

\begin{abstract}
It is proved that for infinitely many positive integers $n$, there exists a circulant graph of order~$n$ whose Weisfeiler-Leman dimension is at least $c\sqrt{\log n}$ for some constant $c>0$ not depending on~$n$.

\end{abstract}

\maketitle

\section{Introduction}\label{120225g}
Under a circulant graph or a \emph{circulant}, we mean a Cayley graph of a cyclic group.
The history of the isomorphism problem for circulants  goes back to the case of a group of prime order, studied in several papers in connection with the \'Ad\'am conjecture (see~\cite{MuzKP2001}). An effective criterion for testing isomorphism suggested by this conjecture was later proved in~\cite{1997a} for cyclic groups of order $n$ not divisible by~$8$ or an odd square. Namely, for any two subsets $X,Y\subseteq C_n$,
\qtnl{280924a}
\cay(C_n,X)\cong \cay(C_n,Y)\quad\Leftrightarrow\quad Y\in X^{\aut(C_n)}
\eqtn
where $X^{\aut(C_n)}=\{X^\sigma:\ \sigma\in\aut(C_n)\}$. The proof of this result and the results mentioned below are based on the method of Schur rings, first applied to the isomorphism problem for circulants in~\cite{KlinP1978}.

For an arbitrary positive integer $n$, the criterion \eqref{280924a} is not true. In two papers~\cite{EvdP2004a} and~\cite{Muz2004}, the  isomorphism problem for circulants was completely and independently solved by constructing polynomial-time  (in~$n$) algorithms. Although the technique used in these two papers was different, the common ground was the use of permutation groups. A natural generalization of this approach was a polynomial-time algorithm for testing isomorphism of circulant objects in an arbitrary concrete category~\cite{MuzP2020}.

In the mid-2000s, at a conference, one of the authors of this paper was asked to briefly describe the idea of an algorithm for testing isomorphism of circulants. This seemed like a difficult task, since the algorithms used (to varying degrees) both the theory of permutation groups and the theory of Schur rings. It was then that the question arose about the existence of a ``purely combinatorial'' polynomial-time algorithm that does not rely on algebraic methods. As a prototype of such an algorithm, it was natural to consider the multidimensional Weisfeiler-Leman algorithm, see~\cite{Kiefer2020}.

For a positive integer $m$ and an arbitrary graph $\cG$ (directed or undirected), the $m$-dimensional Weisfeiler-Leman algorithm (the $m$-dim $\WL$ algorithm) canonically constructs a special coloring of the $m$-tuples consisting of vertices of~$\cG$. Two given graphs are declared isomorphic if the resulting colorings do not distinguish these graphs; in this case, the graphs are said to be \emph{$\WL_m$-equivalent}. The algorithm works correctly if the $\WL_m$-equivalence of the input graphs implies that they are isomorphic. The $m$-dim $\WL$ algorithm is the naive canonical vertex classification in the sense of~\cite{Babai1979a} if  $m=1$,  and is the classical Weisfeiler-Leman algorithm \cite{WLe68} if $m=2$. As one of the key ingredient, the $m$-dim $\WL$ algorithm was used in the Babai quasipolynomial algorithm testing isomorphism of graphs~\cite{Babai2019}.

Following \cite{Grohe2017}, the Weisfeiler-Leman dimension (\emph{$\WL$-dimension}) $\dim_{\scriptscriptstyle\WL}(\cG)$ of a graph $\cG$ is defined to be the smallest positive integer $m$ such that any graph which is $\WL_m$-equivalent to the graph $\cG$ is isomorphic to $\cG$. The isomorphism of any graph of order $n$ and $\WL$-dimension $m$ with any other graph can be tested in time $O(n^{m+1}\log n)$. Now the question we started with can be reformulated as follows (see~\cite{Ponomarenko2022a}): whether there exists a positive integer $M$ such that  $\dim_{\scriptscriptstyle\WL}(\cG)\le M$ for every circulant graph~$\cG$.

It is known that the $\WL$-dimension of a general graph cannot be bounded from above by an absolute constant~\cite{CaiFI1992}. On the other hand, it was proved in~\cite{Ponomarenko2022a} that if $\cG$ is a circulant graph of prime power order, then $\dim_{\scriptscriptstyle\WL}(\cG)\le 3$. The main result of the present paper (Theorem~\ref{090524d} below) shows that the $\WL$-dimension of a general circulant graph can be arbitrary large. To formulate this result, we refer to Section~\ref{180524g} for the exact definition of $\varepsilon$-expanders, and  given a positive integer~$a$, denote by $\fB_a$ the set of all  bipartite cubic $3$-connected graphs of order~$a$.

\thrml{090524d}
Let $\varepsilon>0$ be a real number, and let $a=2b$ be a positive integer. Assume that there is an  $\varepsilon$-expander belonging to $\fB_a$, and let $n=n_1\cdots n_a$, where each~$n_i$ is a product of two distinct primes $\ge 5$ and the products $n_1\cdots n_b$ and $n_{b+1}\cdots n_{2b}$ are coprime. Then there exists a circulant graph of order~$n$ the $\WL$-dimension of which is at least $\lfloor c\cdot \sqrt{a}\rfloor$, where $c=c(\varepsilon)$ is a constant not depending on~$a$, $b$, and~$n$.
\ethrm

Graphs satisfying the conditions of Theorem~\ref{090524d} depend on the choice of the even positive integer~$a$, an $\varepsilon$-expander belonging to $\fB_a$, and the numbers~$n_1,\ldots, n_a$. It is known that for any sufficiently large numbers $a$ and a real number $\varepsilon$, there exist cubic $\varepsilon$-expanders~\cite{Alon2021}. However, we need them to be bipartite and $3$-connected. The existence of such graphs can be deduced using the construction of cubic Ramanujan graphs in~\cite{Morgenstern1994} (see the proof of Theorem~\ref{090524d}). The restrictions on the choice of numbers are quite weak; in particular, they can be chosen so that the number $n$ is square free. Thus the algorithm testing isomorphism of circulants on the base of the criterion~\eqref{280924a} can be simulated by the $m$-dim $\WL$ algorithm for no constant~$m$ independent of~$n$.

\crllrl{250924a}
For infinitely many positive integers $n$, there exists a circulant graph of order~$n$ whose $\WL$-dimension is at least $c\sqrt{\log n}$ for some constant $c>0$ not depending on~$n$.
\ecrllr

Following \cite{Wu2024}, we define the $\WL$-dimension of a circulant $\cG$  with respect to the class of all circulant graphs as  the smallest positive integer $m$ such that any circulant graph which is $\WL_m$-equivalent to~$\cG$ is isomorphic to $\cG$.  It was proved in~\cite{Wu2024} that this ``relative'' $\WL$-dimension of any circulant of order $n$ is at most~$O(\log n)$. The question of the tight asymptotic for the maximum relative $\WL$-dimension of circulants remains open. As Corollary~\ref{250924a} shows, it is between $O(\sqrt{\log n})$ and~$O(\log n)$.

We conclude the introduction with a brief discussion of the proof of Theorem~\ref{090524d}.  Our technique is based on theory of coherent configurations~\cite{CP2019}, see also Section~\ref{091024a}. If $\cG$ is a graph, then the output $\cX=\WL(\cG)$ of the $2$-dim $\WL$-algorithm applied to~$\cG$ can naturally be viewed as a coherent configuration. The concept of the $\WL$-dimension  extends to coherent configurations such that
$$
\dim_{\scriptscriptstyle\WL}(\cG)=\dim_{\scriptscriptstyle\WL}(\cX),
$$
see Section~\ref{040924e}. When the graph $\cG$ is circulant, the coherent configuration~$\cX$ naturally corresponds to a Schur ring  $\cA$ over a cyclic group (a circulant S-ring), see Section~\ref{091024b}. Now to find a circulant graph $\cG$ with large $\dim_{\scriptscriptstyle\WL}(\cG)$, we construct in Section~\ref{190524h} a circulant S-ring $\cA^\star$ such that the corresponding coherent configuration~$\cX^\star$ has two properties: (a) $\dim_{\scriptscriptstyle\WL}(\cX^\star)$ is large, and (b) $\cX^\star=\WL(\cG)$ for some circulant graph~$\cG$. The construction of $\cA^\star$ (see below) enables us to verify property~(b) more or less easily, see  Subsection~\ref{091024d}.

The construction of the S-ring $\cA^*$ satisfying condition~(a) is more subtle. First, we use a theory of multipliers in circulant S-rings, that was developed in~\cite{Evdokimov2016,Evdokimov2015}, see also Section~\ref{101024a}. This theory is applied to express the $\WL_m$-equivalence of coherent configurations of circulant graphs in terms of \emph{$k$-local systems of multipliers} of the corresponding circulant S-rings for $m<\sqrt{k}$, introduced and studied in Section~\ref{101024b}.  Thus our problem  reduces to finding a circulant S-ring having a $k$-local system of multipliers for enough large~$k$. An explicit construction of such an S-ring is given in Section~\ref{101024d} and is based on two ingredients defined in Section~\ref{190524h}.

The first ingredient is a class of coherent configurations, introduced and studied in~\cite{EvdP1999c}, see Subsection~\ref{121123u}.  For any cubic graph $\cG$, one can define a coherent configuration $\cX=\cX(\cG)$ of this class so that: if  $\cG$ has $a$  vertices, then $\cX$ has $4a$ points, and if $\cG$ is a $3$-connected $\varepsilon$-expander, then \qtnl{111024a}
\dim_{\scriptscriptstyle\WL}(\cX)\ge c'\cdot a
\eqtn
for some constant $c'$ depending only on~$\varepsilon$.

The second  ingredient is a class consisting of  special (circulant) S-rings $\cA$ over cyclic groups of orders~$n$ satisfying the conditions of Theorem~\ref{090524d}.  The structure of $\cA$ depends only on the number~$n$, and is enough simple to describe all possible multipliers. In particular, if we fix a bipartite graph $\cG$ from the previous paragraph, then isomorphisms of the coherent configuration $\cX$ yield multipliers of~$\cA$. Now the required S-ring~$\cA^\star$ is the algebraic fusion of the S-ring~$\cA$ with respects to the multipliers corresponding the automorphisms of~$\cX$ (Subsection~\ref{070624a}). Finally, the inequality~\eqref{111024a} enables us to find a $k$-local system of multipliers of~$\cA^\star$ for $k=\lfloor c\cdot a\rfloor$ for some constant $c=c(\varepsilon)$.

\section{Coherent configurations}\label{091024a}
In presenting coherent configurations, we follow the monograph~\cite{CP2019}, which can be consulted for unexplained details.

\subsection{Notation}\label{060225a}
Throughout the paper, $\Omega$ denotes a finite set. For a set $\Delta\subseteq \Omega$, the Cartesian product $\Delta\times\Delta$ and its diagonal are denoted by~$\bone_\Delta$ and $1_\Delta$, respectively. If $\Delta=\{\alpha\}$, we abbreviate $1_\alpha:=1_{\{\alpha\}}$. For a relation $s\subseteq\bone_\Omega$, we set $s^*=\{(\alpha,\beta): (\beta,\alpha)\in s\}$, $\alpha s=\{\beta\in\Omega:\ (\alpha,\beta)\in s\}$ for all $\alpha\in\Omega$, and define $\grp{s}$ as the minimal (with respect to inclusion) equivalence relation on~$\Omega$, containing~$s$. The \emph{support} of~$s$ is defined to be the smallest set $\Delta\subseteq\Omega$ such that $s\subseteq\bone_\Delta$.

For any collection $S$ of relations, we denote by $S^\cup$ the set of all unions of elements of~$S$, and consider $S^\cup$ as a poset with respect to inclusion. The \emph{dot product} of relations $r,s\subseteq\Omega\times\Omega$ is defined by the formula $$r\cdot s:=\{(\alpha,\beta)\!:\ (\alpha,\gamma)\in r,\ (\gamma,\beta)\in s\text{ for some }\gamma\in\Omega\}.$$

The set of classes of  an equivalence relation $e$ on~$\Omega$ is denoted by $\Omega/e$. For $\Delta\subseteq\Omega$, we set $\Delta/e=\Delta/e_\Delta$, where $e_\Delta=\bone_\Delta\cap e$. If the classes of $e_\Delta$ are singletons, $\Delta/e$ is identified with $\Delta$. Given a relation $s\subseteq \bone_\Omega$, we put
\qtnl{130622a}
s_{\Delta/e}=\{(\Gamma,\Gamma')\in \bone_{\Delta/e}:\  s^{}_{\Gamma,\Gamma'}\ne\varnothing\},
\eqtn
where $s^{}_{\Gamma,\Gamma'}=s\cap (\Gamma\times \Gamma')$, and  abbreviate $s^{}_\Gamma:=s^{}_{\Gamma,\Gamma}$. For a fixed $s$, the set of all  equivalence relations $e$ on~$\Omega$, such that $e\cdot s=s\cdot e=s$, contains  the largest (with respect to inclusion) element which is denoted by  $\rad(s)$ and called the {\it radical} of~$s$. Obviously, $\rad(s)\subseteq\grp{s}$.

For a set $B$ of bijections $f:\Omega\to\Omega'$, subsets $\Delta\subseteq \Omega$ and $\Delta'\subseteq \Omega'$, equivalence relations $e$ and $e'$ on $\Omega$ and  $\Omega'$, respectively, we put
$$
B^{\Delta/e,\Delta'/e'}=\{f^{\Delta/e}:\ f\in B,\ \Delta^f=\Delta',\ e^f=e'\},
$$
where $f^{\Delta/e}$ is the bijection from $\Delta/e$ onto $\Delta'/e'$, induced by~$f$; we also  abbreviate $B^{\Delta/e}:=B^{\Delta/e,\Delta'/e}$ if $\Delta'$ is clear from context.

Under a tuple on $\Omega$, we mean an $m$-tuple $x=(x_1,\ldots,x_m)$ for some~$m\ge 1$. We put  $\Omega(x):=\{x_1,\ldots,x_m\}$. For any $1\le k\le m$, we put $\pr_k(x)=(x_1,\ldots,x_k)$ and extend this notation to $\pr_k  \sX$ for all $\sX\subseteq\Omega^m$.

For a finite group $G$ and a set $X\subseteq G$, we treat $\und{X}=\sum_{x\in X}x$ as an element of the group ring $\mZ G$. We put $\rad(X)=\{g\in G:\ gX=Xg=X\}$, denote by~$\grp{X}$  the subgroup of $G$ generated by~$X$, and set $X^{-1}=\{x^{-1}:\ x\in X\}$.

\subsection{Coherent configurations}\label{270224b}
Let $S$ be a partition of $\bone_\Omega$. A pair $\mathcal{X}=(\Omega,S)$ is called a \emph{coherent configuration} on $\Omega$ if
\nmrt
\tm{CC1}  $1_\Omega\in S^\cup$,
\tm{CC2} $s^*\in S$ for all $s\in S$,
\tm{CC3} given $r,s,t\in S$, the number $c_{rs}^t=|\alpha r\cap \beta s^{*}|$ does not depend on $(\alpha,\beta)\in t$.
\enmrt
The number $|\Omega|$ is called the {\it degree} of~$\cX$. A coherent configuration~$\cX$ is said to be {\it trivial} if  $S=S(\cX)$ consists of~$1_\Omega$ and its complement (unless $\Omega$ is a singleton), \emph{discrete} if~$S$ consists of  singletons, and {\it homogeneous} or a {\it scheme} if $1_\Omega\in S$.

With any permutation group $G\le\sym(\Omega)$, one can associate a coherent configuration $\inv(G)=(\Omega,S)$, where $S$ consists of the relations $(\alpha,\beta)^G=\{(\alpha^g,\beta^g): \ g\in G\}$ for all $\alpha,\beta\in\Omega$. A scheme $\cX$ is said to be \emph{regular} if $\cX=\inv(G)$, where the group $G$ is regular, or equivalently, if $|\alpha s|=1$ for all $\alpha\in\Omega$ and $s\in S$.

\subsection{Isomorphisms and algebraic isomorphisms}\label{110622w} A  {\it combinatorial isomorph\-ism} or, briefly, \emph{isomorphism} from a coherent configuration $\cX=(\Omega, S)$ to a coherent configuration $\cX'=(\Omega', S')$ is defined to be a bijection $f: \Omega\rightarrow \Omega'$ such that  $S^f=S'$. In this case,  $\cX$ and $\cX'$ are said to be {\it isomorphic}; the set of all isomorphisms from~$\cX$ to~$\cX'$ is denoted by $\iso(\cX,\cX')$, and by $\iso(\cX)$ if $\cX=\cX'$. The group $\iso(\cX)$ contains a normal subgroup
$$
\aut(\cX)=\{f\in\sym(\Omega):\ s^f=s \text{ for all } s\in S\}
$$
called the {\it automorphism group} of the coherent configuration~$\cX$.

A bijection $\varphi:S\rightarrow S'$ is called an \emph{algebraic isomorphism} from $\cX$ to $\cX'$ if for all $r,s,t\in S$, we have
$$
c_{\varphi (r)\varphi (s)}^{\varphi(t)}=c_{rs}^t.
$$
In this case, $|\Omega'|=|\Omega|$ and $|S|=|S'|$. The set of all algebraic isomorphisms from~$\cX$ to $\cX'$ is denoted by $\iso_{alg}(\cX,\cX')$, and by $\aut_{alg}(\cX)$ if $\cX=\cX'$.

Every combinatorial  isomorphism $f\in\iso(\cX,\cX')$ induces the algebraic isomorphism $\varphi_f:\cX\to \cX',$ $s\mapsto s^f$; we put
$$
\iso(\cX,\cX',\varphi)=\{f\in\iso(\cX,\cX'):\ \varphi_f=\varphi\}
$$
and abbreviate $\iso(\cX,\varphi):=\iso(\cX,\cX,\varphi)$. Note that $\aut(\cX)=\iso(\cX,\id_\cX)$, where $\id_\cX$ is the trivial (identical) algebraic automorphism of~$\cX$. If  $\cX$  is a regular scheme, then $\iso(\cX,\cX',\varphi)\ne\varnothing$ for any~$\varphi$ (see~\cite[Theorem~2.3.33]{CP2019}).

\subsection{Relations and fibers}  Let $\cX=(\Omega,S)$ be a coherent configuration. The elements of $S$ and of $S^\cup$ are called {\it basis relations} and \emph{relations} of  ~$\cX$. The basis relation containing a pair $(\alpha,\beta)\in\bone_\Omega$ is denoted by $r_\cX(\alpha,\beta)$ or simply by $r(\alpha,\beta)$. The set of all relations of $\cX$ is closed with respect to intersections, unions, and the dot product.

Any set $\Delta\subseteq\Omega$ such that $1_\Delta\in S$ is called a {\it fiber} of $\cX$. In view of condition~(CC1), the set~$F=F(\cX)$  of all fibers forms a partition of~$\Omega$. Every basis relation is contained in the Cartesian product of two uniquely determined fibers. 

Every algebraic isomorphism $\varphi\in\iso_{alg}(\cX,\cX')$ can be extended in a natural way to a poset isomorphism $S^\cup\to (S')^\cup$, denoted also by~$\varphi$; we have $1^{}_{\Omega'}=\varphi(1_\Omega)$. This poset isomorphism respects the set theoretical operations and the dot product. Furthermore, $\varphi$ induces a  poset isomorphism $F^\cup\to (F')^\cup$ which preserves fibers. It takes  a set $\Delta$ to the set $\Delta^\varphi$, and is defined by the equality $\varphi(1_\Delta)=1_{\Delta^\varphi}$.


\subsection{Parabolics.}\label{ssect:parabolics}
An equivalence relation $e$ on a set $\Delta\subseteq\Omega$ is called a \emph{partial parabolic} of the coherent configuration~$\cX$ if $e$ is a relation of~$\cX$; if, in addition, $\Delta=\Omega$, then $e$ is called a \emph{parabolic} of~$\cX$. Note that the reflexive and  transitive closure of any symmetric relation of $\cX$ is a partial parabolic.


Let $e$ be a partial parabolic on~$\Omega$, and let $\Delta\in\Omega/e$. Denote by $S_{\Omega/e}$ and $S_\Delta$ the sets of all non-empty relations
$$
s_{\Omega/e}=\{(\Delta,\Gamma)\in\Omega/e\times\Omega/e:\ s_{\Delta,\Gamma}\ne\varnothing\}\quad\text{and}\quad s_{\Delta}=s_{\Delta,\Delta},
$$
respectively, where $s\in S$. Then the pairs $\cX_{\Omega/e}=(\Omega/e,S_{\Omega/e})$ and $\cX_\Delta=(\Delta,S_\Delta)$ are coherent configurations called the \emph{quotient} of~$\cX$ modulo~$e$ and \emph{restriction} of~$\cX$ to~$\Delta$. Note that if $e=1_\Delta$ for some $\Delta\in F^\cup$, then  $\cX_{\Omega/e}=\cX_\Delta$.

Let $\varphi\in\iso_{alg}(\cX,\cX')$ be an algebraic isomorphism and  $e$ a (partial) parabolic of~$\cX$. Then  $\varphi(e)$ is a (partial) parabolic of~$\cX'$ with the same number of classes. Moreover,
\qtnl{030322a}
\varphi(\grp{s})=\grp{\varphi(s)}\qaq \varphi(\rad(s))=\rad(\varphi(s)).
\eqtn
The algebraic isomorphism $\varphi$ induces the algebraic isomorphism
\qtnl{060524a}
\varphi^{}_{\Omega^{}/e^{}}:\cX^{}_{\Omega^{}/e^{}}\to \cX'_{\Omega'/e'},\ s^{}_{\Omega^{}/e^{}}\mapsto s'_{\Omega'/e'},
\eqtn
where $e'=\varphi(e)$ and $s'=\varphi(s)$.  Now let  $\Delta\in\Omega/e$ and $\Delta'\in\Omega/e'$. Assume that there is $\Gamma\in F$ such that $\Delta\cap\Gamma\ne\varnothing\ne \Delta'\cap\Gamma'$, where $\Gamma'=\Gamma^\varphi$. Then  $\varphi$ induces an algebraic isomorphism
\qtnl{200322u}
\varphi^{}_{\Delta,\Delta'}:\cX^{}_{\Delta^{}}\to \cX'_{\Delta'},\ s^{}_{\Delta^{}}\mapsto s'_{\Delta'}.
\eqtn
This  algebraic isomorphism always exists for all $\Delta$ and $\Delta'$ if $\cX$ (and hence $\cX'$) is a scheme.

\subsection{Tensor and wreath product.}
Let $\cX=(\Omega,S)$ and $\cX'=(\Omega',S')$ be two coherent configurations. Denote by $S\otimes S'$ the set of all relations
$$
s\otimes s'=\left\{((\alpha,\alpha'),(\beta,\beta'))\in (\Omega\times\Omega')^2:\ (\alpha,\beta)\in s,\ (\alpha',\beta')\in s'\right\},
$$
where $s\in S$ and  $s'\in S'$. Then the pair $\cX\otimes\cX'=(\Omega\times \Omega',S\otimes S')$ is a coherent configuration called the \emph{tensor product} of~$\cX$ and~$\cX'$.  One can easily prove that $\aut(\cX\times\cX')=\aut(\cX)\times\aut(\cX')$.

Assume that the coherent configurations $\cX$ and $\cX'$ are schemes. Denote by $T$  and $T'$ the sets of all relations
$$
s\otimes 1_{\Omega'},\ s\in S, \quad\qaq\quad\bone_\Omega\otimes s',\ 1_{\Omega'}\ne s'\in S',
$$
respectively. Then  the pair $\cX\wr\cX'=(\Omega\times \Omega',T\,\cup\, T')$ is also a scheme. It is called the \emph{wreath  product} of~$\cX$ and~$\cX'$. The mapping $\aut(\cX\wr\cX')\to\aut(\cX')$ , $f\mapsto f^{\Omega'},$ is an epimorphism the kernel of which is equal to the direct product of $|\Omega'|$ copies of the group $\aut(\cX)$.

\section{Binary  schemes}
In the proof of our main result, we will need schemes for which every ``correct'' bijection between point tuples of bounded arity lifts to an automorphism. The automorphism groups of these schemes are the binary permutation groups in sense of~\cite{Cherlin2000}.

Let $\cX$ be a scheme on $\Omega$. For any integer $m\ge 1$ and a tuple $x\in \Omega^m$, we define the $m\times m$ array
\qtnl{170423a}
R_{\cX,m}(x)=(r_\cX(x_i,x_j):\ 1\le i,j\le m).
\eqtn
The scheme $\cX$ is said to be \emph{binary} if for any $m$ and any two $m$-tuples $x,y\in \Omega^m$ such that $R_{\cX,m}(x)=R_{\cX,m}(y)$, there exists an automorphism of $\cX$, that takes $x$ to~$y$. In this case, the group $\aut(\cX)$ is transitive, and acts transitively on each basis relation of~$\cX$.

The main goal of this section is to prove that any scheme constructed from regular schemes using tensor and wreath products is binary.

\thrml{250724a}
The class of binary schemes contains all regular schemes and is closed with respect to tensor and wreath products.
\ethrm
\prf
Below, $\cX$ is a scheme on $\Omega$, $m\ge 1$ is an integer, and $x,y\in \Omega^m$ are such that $R(x)=R(y)$, where $R=R_{\cX,m}$. First, assume that $\cX$ is regular. Then there exists a unique 
automorphism $f\in\aut(\cX)$ such that ${x_1}^f=y_1$. Let $1\le i\le m$. By the assumption, the point  $y_i$ is the unique one  for which $r(x_1,x_i)=r(y_1,y_i)$, where $r=r_\cX$. Since $r(x_1,x_i)^f =r(y_1,x_i^f)$, this shows that $y^{}_i={x_i}^f$. Thus $y=x^f$ and $\cX$  is binary.

Let $\cX=\cX_1\otimes\cX_2$, where $\cX_1$ and $\cX_2$ are  binary schemes on~$\Omega_1$ and $\Omega_2$, respectively. Denote by $\pr_i:\Omega_1\times\Omega_2\to\Omega_i$ the natural projection, $i=1,2$.
Identifying  in a natural way the sets $(\Omega_1\times\Omega_2)^m$ and $\Omega_1^{\phantom{\,}m}\times\Omega_2^{\phantom{\,}m}$, we can write   $x=(\pr_1 x,\pr_2x)$ and  $y=(\pr_1 y,\pr_2y)$.
 By the definition of tensor product, we have
$$
R_1(\pr_1 x)=R_1(\pr_1y)\qaq R_2(\pr_2 x)=R_2(\pr_2y),
$$
where $R_1=R_{\cX_1,m}$,  and $R_2=R_{\cX_2,m}$. Since the scheme $\cX_i$ is binary,  there exists an automorphism $f_i\in\aut(\cX_i)$ such that $\pr_iy=(\pr_i x)^{f_i}$, $i=1,2$.  It follows that the permutation $f:=(f_1,f_2)$ of~$\Omega_1\times\Omega_2$ belongs to the group $\aut(\cX_1)\times\aut(\cX_2)=\aut(\cX)$. Thus,
$$
x^f=(\pr_1 x,\pr_2 x)^{(f_1,f_2)}=((\pr_1 x)^{f_1},(\pr_2 x)^{f_2})=(\pr_1y,\pr_2y)=y,
$$
and the scheme $\cX$  is binary.

Now let $\cX=\cX_1\wr\cX_2$. Then $e=\bone_{\Omega_1}\otimes 1_{\Omega_2}$ is a parabolic of the scheme $\cX$. Each class of~$e$ is of the form $\bar\alpha=\Omega_1\times\{\alpha\}$, where $\alpha\in\Omega_2$. The  bijection $\Omega_2\to \Omega/e$, $\alpha\mapsto\bar\alpha$, $\alpha\in\Omega_2$ is an isomorphism from the scheme $\cX_2$ to the quotient scheme $\bar\cX=\cX_{\Omega/e}$; in particular, the scheme $\bar \cX$ is binary.
$$
\bar R(\bar x)=\bar R(\bar y),
$$
where $\bar R=R_{\bar\cX,m}$, $\bar x=(\ov{ x_1},\ldots,\ov{x_m})$, and  $\bar y=(\ov{ y_1},\ldots,\ov{y_m})$. Thus  there exists an automorphism $\bar f\in\aut(\bar \cX)$ such that $\bar y={\bar x}^{\bar f}$. Because the natural homomorphism $\aut(\cX)\to\aut(\bar\cX)$ is surjective (see~\cite[Theorem~3.4.6.]{CP2019}), one can find an automorphism $g\in\aut(\cX)$ such that $g^{\Omega/e}=\bar g$. Replacing $x$ by $x^g$ and taking into account that
$R(x^g)=R(x)$, we may assume that
\qtnl{210824c}
\ov{x_i}=\ov{y_i},\qquad i=1,\ldots,m.
\eqtn

Let $\Delta\in \Omega/e$. Let the integer $k\ge 1$ and indices  $1\le i_1<\cdots <i_k\le m$ be such that $\{x_{i_1},\ldots x_{i_k}\}=\Omega(x)\cap\Delta$. Put  $x_\Delta=(x_{i_1},\ldots, x_{i_k})$. In a similar way, we define the integer $\ell\ge 1$ and $\ell$-tuple~$y_\Delta$. Since $R(x)=R(y)$, formula~\eqref{210824c} implies that $k=\ell$ and
$$
R_\Delta(x_\Delta)=R_\Delta(y_\Delta),
$$
where $R_\Delta=R_{\cX_\Delta,k}$. Because the scheme $\cX_\Delta=\cX_1$ is binary, there exists an automorphism $f_\Delta\in\aut(\cX_\Delta)$ such that $y_\Delta=(x_\Delta)^{f_\Delta}$. Let us define a permutation $f\in\sym(\Omega)$ such that $f^\Delta=f_\Delta$ for all $\Delta\in\Omega/e$. Then obviously $x^f=y$. To complete the proof, it remains to note that $f\in\aut(\cX)$, because the kernel of the homomorphism $\aut(\cX)\to\aut(\ov\cX)$  coincides with the direct product of $|\Omega/e|$ copies of the group $\aut(\cX_\Delta)$.
\eprf

\section{WL-dimension of graphs and coherent configurations}\label{040924e}
\subsection{Definitions} Let $\sR$ be a linear ordered set of binary relations on $\Omega$, and let $m\ge 1$ be an integer. The Weisfeiler-Leman algorithm (the $m$-dim $\WL$) constructs  a coloring $c=c_{m,\sR}$ assigning a color $c(x)$ to each $m$-tuple $x\in\Omega^m$. The coloring is canonical in the sense that for any bijection~$f:\Omega\to\Omega'$, the induced bijection $\Omega^m\to {\Omega'}\phmaa{m}$ takes the coloring $c_{m,\sR}$ to the coloring  $c_{m,\sR'}$, where $\sR'=\{s^f:\ s\in\sR\}$. More details (including the exact definition of the $m$-dim $\WL$) can be found in~\cite{Grohe2017} (see also~\cite{AndresHelfgott2017}).

The  color classes of the coloring $c_{m,\sR}$ form a partition  $\WL_m(\sR)$ of the Cartesian power~$\Omega^m$,\footnote{Though the coloring $c_{m,\sR}$ depends on the linear order of~$\sR$, the partition $\WL_m(\sR)$ does not.}  which is an  $m$-ary coherent configuration in the sense of~\cite{Babai2019}; for $m=2$, the $m$-ary coherent configurations are just the coherent configurations  in the sense of Subsection~\ref{270224b}. We do not need to give here the exact definition of  $m$-ary coherent configurations or their algebraic isomorphisms; the interested reader can find them in~\cite[Section~3]{Ponomarenko2022a}. Let us make just a few remarks about projections.

Let $\fX$ be an $m$-ary coherent configuration and $1\le k\le m$.  The collection of the projections $\pr_k\sX$, $\sX\in \fX$, forms a partition of $\Omega^k$, that is a $k$-ary coherent configuration; we denote it by $\pr_k \fX$. The projections are respected by algebraic isomorphisms. More precisely, any algebraic isomorphism  $\wh\varphi:\fX\to\fX'$ of $m$-ary coherent configurations induces the algebraic isomorphism $\wh\varphi_k:\pr_k \fX\to\pr_k \fX'$ of $k$-ary coherent configurations such that $\pr_k \wh\varphi(\sX)=\wh\varphi_k(\pr_k\sX)$ for all $\sX\in\fX$.

At this point, we assume that $m\ge 2$.  For a ($2$-ary) coherent configuration $\cX$, denote by $\WL_m(\cX)$ the $m$-ary coherent configuration $\WL_m(\sR)$ with $\sR=S(\cX)$. Following~\cite{Ponomarenko2022a}, we say that two coherent configurations $\cX$ and $\cX'$ are   \emph{$\WL_m$-equivalent} with respect to an algebraic isomorphism~$\varphi:\cX\to \cX'$  if there exists an algebraic isomorphism  $\hat\varphi:\WL_m(\cX)\to\WL_m(\cX')$  such that its $2$-projection $\hat\varphi_2$ extends~$\varphi$, i.e., $\hat\varphi_2(s)=\varphi(s)$ for all $s\in S(\cX)$.
(Note that we always have $\pr_2 \WL_m(\cX)\ge \cX$ and $\pr_2 \WL_m(\cX')\ge \cX'$.)

The $m$-dim $\WL$ algorithm \emph{identifies} a coherent configuration~$\cX$ if  for any algebraic isomorphism  $\varphi:\cX\to\cX'$ with respect to which  $\cX$ and $\cX'$ are $\WL_m$-equivalent, the set $\iso(\cX,\cX',\varphi)$ is not empty. The \emph{$\WL$-dimension} $\dim_{\scriptscriptstyle\WL}(\cX)$ of the coherent configuration~$\cX$ is defined to be the smallest positive integer $m\ge 2$ such that the $m$-dim~$\WL$ algorithm identifies~$\cX$.

With any graph $\cG$, we associate the coherent configuration $\WL(\cG)=\WL_2(\sR)$, where the set $\sR$ is the singleton consisting of the edge set of $\cG$. Let us define the $\WL$-dimension of~$\cG$ by the formula
\qtnl{040924a}
\dim_{\scriptscriptstyle\WL}(\cG)=\dim_{\scriptscriptstyle\WL}(\WL(\cG)).
\eqtn
It should be noted that our definition of the $\WL$-dimension of a graph differs somewhat from the generally accepted one.  However, both definitions give the same value of the $\WL$-dimension when the graph in question is of $\WL$-dimension at least~$2$ in the conventional sense. Therefore, the difference in the definitions of the $\WL$-dimension does not affect the validity of the main result of our paper.

\subsection{Pebble games}
It is convenient to express the $\WL_m$-equivalence of coherent  configurations $\cX$ and $\cX'$ with respect to an algebraic isomorphism $\varphi:\cX\to \cX'$ in terms of the  pebbling game $\sC_m(\varphi)$, see~\cite{Wu2024,Chen2023a}.  The idea and the main result of this subsection go back to a fundamental paper~\cite{CaiFI1992}.

There are two players, called Spoiler and Duplicator, and $m$ pairwise distinct pebbles, each given in duplicate.   The game consists of rounds and each round consists of two parts. At the first part, Spoiler chooses a set~$A'$ of points belonging to~$\cX$ or~$\cX'$. Duplicator responds with a set $A$ in the other coherent configuration, such that $|A|=|A'|$ (if this is impossible, then Duplicator loses). At the second part, Spoiler places one of the pebbles\footnote{It is allowed to move previously placed pebbles to other 	vertices and place more than one pebble on the same vertex.}  on a point in~$A$. Duplicator responds  by placing the copy of the pebble on some point of~$A'$.

The configuration  after a round is determined by a bijection $f:\Delta\to\Delta'$, where $\Delta\subseteq\Omega$ (respectively, $\Delta'\subseteq\Omega' $) is the set of points  in $\cX$ (respectively, in~$\cX'$), covered by pebbles, and any two points $\alpha\in\Delta$ and $\alpha^f\in\Delta'$ are covered by the copies of the same pebble. Duplicator wins the round if
$$
(s_{\Delta^{}})^f=\varphi(s)_{\Delta'}
$$
for all $s\in S$ such that $s_\Delta\ne\varnothing$. Spoiler wins if Duplicator  does not win.

The game starts from an initial configuration (which  is considered as the configuration after zero rounds), i.e.,  a pair $(x,x')\in\Omega^k\times {\Omega'}\phmaa{k}$  for which $k\le m$ and the points~$x^{}_i$ and~$x'_i$ are covered with copies of the same pebble, $i=1,\ldots,k$. We say that Duplicator  has a winning strategy  for the game $\sC_m(\varphi)$  on $\cX$ and $\cX'$ with  initial configuration $(x,x')$ if, regardless of Spoiler's actions, Duplicator wins after any number of rounds.

\lmml{040924b}
Let $\cX$ and $\cX'$ be coherent configurations on, respectively, $\Omega$ and $\Omega'$, and let $\varphi:\cX\to\cX'$ be an algebraic isomorphism. Assume that  an integer $m\ge 2$ and a bijection $f:\Omega^m\to {\Omega'}\phmaa{m}$ are such that for each  $x\in\Omega^m$,  Duplicator has a winning strategy in the  game  $\sC_{m+1}(\varphi)$ on $\cX$ and $\cX'$ with initial configuration $(x,x')$, where $x'=x^f$. Then $\cX$ and $\cX'$ are $\WL_m$-equivalent with respect to~$\varphi$.
\elmm
\prf
Let $T\subseteq \Omega^m\times {\Omega'}\phmaa{m}$ consist  of all pairs $(x,x')$ such that Duplicator has a winning strategy in the pebble game  $\sC_{m+1}(\varphi)$ on $\cX$ and $\cX'$ with initial configuration~$(x,x')$.  By the assumption, $T$ contains each pair $(x,x^f)$, $x\in\Omega^m$. Since $f$ is a bijection,  the relation $T$ has full support. This completes the proof, because, by \cite[Theorem~3.3]{Wu2024},  $\cX$ and $\cX'$ are $\WL_m$-equivalent with respect to~$\varphi$ if and only if $T$ has  full support.
\eprf

\section{S-rings and Cayley schemes}\label{091024b}

\subsection{S-rings} Let $G$ be a finite group. A free $\mZ$-submodule~$\cA$ of the group ring~$\mZ G$ is  called a {\it Schur ring} ({\it S-ring}, for short) over~$G$ if there exists a partition $\cS=\cS(\cA)$ of~$G$ such that $\cA=\Span_\mZ\{\und{X}:\ X\in\cS\}$ and
\nmrt
\tm{S1} $\{1_G\}\in\cS$,
\tm{S2} $X^{-1}\in\cS$ for all $X\in\cS$,
\tm{S3} $\und{X}\,\und{Y}=\sum_{Z\in\cS}c_{X,Y}^Z\und{Z}$ for all $X,Y\in \cS$ and some nonnegative  integers $c_{X,Y}^Z$.
\enmrt
The elements of $\cS$ and the number $\rk(\cA)=|\cS|$ are called, respectively, the {\it basic sets} and {\it rank} of~$\cA$. The  integer $c_{X,Y}^Z$ is equal to the number of representations $z=xy$ with $x\in X$ and $y\in Y$ for a fixed $z\in Z$. For any subgroup $M\le\aut(G)$, the set $\{x^M:\ x\in G\}$ of all $M$-orbits forms a partition  $\cS(\cA)$ for some  S-ring~$\cA$  over~$G$; we denote this S-ring by $\cyc(M,G)$ and call it \emph{cyclotomic}.

The partial order  $\le$ on the  S-rings over $G$ is induced by inclusion. Thus, $\cA\le\cA'$ if and only if any basic set of $\cA$ is a union of some basic sets of $\cA'$; in this case, we say that $\cA'$ is an \emph{extension} of~$\cA$. The least and greatest elements with respect to~$\le$ are, respectively,  the \emph{trivial}  S-ring spanned by~$\und{1_G}$ and~$\und{G}$, and the group ring~$\mZ G$. For any $X\subseteq G$, we denote by $\WL_G(\und{X})$ the smallest S-ring over~$G$, that contains the element~$\und{X}$.

\subsection{Isomorphisms}\label{080522f}
For any set $X\subseteq G$, one can define a Cayley (di)graph $\cay(G,X)$ with vertex set $G$ and the arcs $(x,y) $ with $yx^{-1}\in X$. The automorphism group of this graph contains a regular subgroup $G_{right}$ of the symmetric group $\sym(G)$, induced by right multiplications of~$G$.

Let $\cA$  and $\cA'$ be S-rings over groups $G$ and $G'$, respectively, and let $f:G\to G'$ be a bijection. We say that $f$ is  a {\it (combinatorial) isomorphism} from $\cA$ to $\cA'$ if  for each $X\in\cS(\cA)$, there is $X'\in\cS(\cA')$ such that $f$ is a graph isomorphism from $\cay(G,X)$ to $ \cay(G',X')$, i.e.,
$$
y\,x^{-1}\in X\quad \Leftrightarrow\quad y' \,{x'}^{-1}\in X'
$$
for all $x,y\in G$, where $x'=x^f$ and $y'=y^f$.\footnote{This condition means exactly that $(Xx)^f=X'x'$ for all $X\in\cS$ and $x\in G$.}  The isomorphism~$f$ is said to be \emph{normalized} if $f(1_{G^{}})=1_{G'}$. Note that the S-rings~$\cA$ and~$\cA'$ are isomorphic if and only if there is a normalized isomorphism from~$\cA$ to~$\cA'$.

The group of all isomorphisms from the S-ring~$\cA$ to itself is denoted by $\iso(\cA)$. It has a normal subgroup  $\aut(\cA)$, the \emph{automorphism group} of~$\cA$,  that is equal to the intersection of the groups $\aut(\cay(G,X))$, $X\in\cS$.  In particular,
\qtnl{311222a}
G_{right}\le \aut(\cA)\le\aut(\cay(G,X))
\eqtn
for each~$X$. Thus, $f\in\aut(\cA)$ if and only if for all $x,y\in G$ and all $X\in\cS$, the elements $y\,x^{-1}$ and  $y'\,{x'}^{-1}$ belong or do not belong to~$X$ simultaneously. Note that if $f\in\aut(G)$ or $f$ is a normalized isomorphism of~$\cA$, then $f\in\aut(\cA)$ if and only if $f$ leaves every basic set of~$\cA$ fixed.

\subsection{$\cA$-sets}
Any union of basic sets is called an \emph{$\cA$-set}. The set of all of them is closed with respect to taking inverse,  product, and standard set-theoretical operations. An $\cA$-set which is a subgroup of~$G$ is called an {\it $\cA$-subgroup}. For example, if~$X$ is an $\cA$-set, then  the groups $\grp{X}$ and $\rad(X)$ are $\cA$-subgroups.  

The following important theorem goes back to I.~Schur and H.~Wielandt (see \cite[Ch.~IV]{Wielandt1964}); as to the formulation given here we refer to~\cite[Theorem~2.3]{Evdokimov2013}. Below for a subset $X$ of an abelian  group $G$, an integer $m$, and a prime $p$, we put
$$
X^{(m)}=\{x^m:\ x\in X\},\quad
X^{[p]}=\{x^p:\ x\in X,\ |xH\cap X|\not\equiv 0\mod p\}
$$
where $H=\{g\in G:\ g^p=1\}$.

\thrml{261009b}
Let $\cA$ be an S-ring over an abelian  group of order~$n$, and $\cS=\cS(\cA)$. Then
\nmrt
\tm{1} $X^{(m)}\in\cS$ for any integer $m$ coprime to $n$ and any $X\in\cS$,
\tm{2} $X^{[p]}\in\cS^\cup$ for any prime $p$ dividing $n$ and any $X\in\cS^\cup$.\footnote{Though Lemma~2.1 in \cite{EvdP2005} requires $X$ to be the basic set of $\cA'$, the proof of the lemma remains true under a weaker assumption that $X$ is an $\cA'$-set.}
\enmrt
\ethrm

\subsection{Tensor and wreath products}\label{150225b}
Let $\cA_1$ and $\cA_2$ be S-rings over groups $G_1$ and $G_2$, respectively. The \emph{tensor product} $\cA_1\otimes\cA_2$ and the \emph{wreath product} $\cA_1\wr\cA_2$ are defined to be the S-rings over the direct product $G=G_1\times G_2$, such that
$$
\cS(\cA_1\otimes\cA_2)=\{X_1\times X_2: X_1\in\cS_1,\ X_2\in\cS_2\},
$$
where $\cS_1=\cS(\cA_1)$ and $\cS_2=\cS(\cA_2)$, and
$$
\cS(\cA_1\wr\cA_2)=\{X_1\times \{1_{G_2}\}: X_1\in\cS_1\}\, \cup\, \{G_1\times X_2: 1_{G_2}\ne X_2\in\cS_2\},
$$
respectively. If $H_1\le G_1$ and $H_2\le G_2$ are, respectively, $\cA_1$- and $\cA_2$-subgroups, then $H_1H_2$ is an $(\cA_1\otimes\cA_2)$-subgroup, whereas $H_1$ and $G_1H_2$ are $(\cA_1\wr\cA_2)$-subgroups. If, in addition, the group $G_1\times G_2$ is cyclic, then these are the only $(\cA_1\otimes\cA_2)$- and $(\cA_1\wr\cA_2)$-subgroups, respectively.

\subsection{Sections}
Let $U$ and $L$ be $\cA$-subgroups and also $L$  a normal subgroup of $U$. The group $\fs=U/L$ is called an \emph{$\cA$-section}; the set of all $\cA$-sections is denoted by~$\fS(\cA)$. Furthermore, let $\pi:U\to U/L$ be the natural epimorphism. The set
$$
\cS_\fs=\{\pi(X):\ X\in\cS,\ X\subseteq U\}
$$
forms a partition of  the group $\fs$, and the elements $\und{\pi(X)}$ span an S-ring over $\fs$; it is denoted by $\cA_\fs$. If $U=\grp{X}$ and $L=\rad(X)$ for some $X\in\cS$, then the section $\fs$ is said to be \emph{principal} and is denoted by $\fs_X$.  In this case, setting $\fs=\fs_X$,  $\pi_\fs=\pi$, and  $X_\fs=\pi(X)$, we have
\qtnl{240824h}
X={\pi_\fs}^{-1}(X_\fs).
\eqtn
For any $x\in U$, we put $x_\fs=\pi_\fs(x)$.

Let $\fs=U/L$ and $\ft=U'/L'$ be $\cA$-sections. We say that $\fs$  is  a {\it subsection} of $\ft$  if $U\le U'$ and $L\ge L'$; in this case, we write $\fs\preceq \ft$. This defines a partial order on the set  $\fS(\cA)$. The greatest element in it is $G/1$;  every minimal element is of the form $H/H$, where $H$ is an $\cA$-subgroup.

\subsection{Projective equivalence} An $\cA$-section $\ft=U'/L'$ is called a {\it multiple} of the section~$\fs=U/L$ if
\qtnl{310114b}
L=U\cap L'\quad\text{and}\quad U'=UL'.
\eqtn
The {\it projective equivalence} relation on the set $\fS(\cA)$ of all $\cA$-sections is defined to be the symmetric transitive closure  of the relation ``to be a multiple"; we write $\fs\sim\ft$ if the sections $\fs$ and $\ft$ are projectively equivalent. Any two projectively equivalent sections are obviously isomorphic as groups. Moreover, if  $\fs\sim\ft$, then there is a group isomorphism $f=f_{\fs,\ft}$ from~$\fs$ to $\ft$ belonging to $\iso(\cA_\fs,\cA_\ft)$ and such that $(\gamma^\fs)^f=\gamma^\ft$ for all normalized automorphisms
$\gamma\in\aut(\cA)$, where $\gamma^\fs$ and $\gamma^\ft$ are the restriction of $\gamma$ to~$\fs$ and to~$\ft$, respectively (see Subsection~\ref{060225a}).

\subsection{Algebraic isomorphisms}
Let $\cA$ and $\cA'$ be S-rings over groups $G$ and $G'$, respectively. A ring isomorphism $\varphi:\cA\to\cA'$ is called an \emph{\it algebraic isomorphism} from $\cA$ to $\cA'$, if for any $X\in\cS(\cA)$, there exists $X'\in\cS(\cA')$ such that
$$
\varphi(\und{X})=\und{X'}.
$$
Each normalized isomorphism $f:\cA\to\cA'$ \emph{induces} an algebraic isomorphism  for which $X'=X^f$. However not every algebraic isomorphism is induced by an isomorphism.

From the definition, it follows that the mapping $X\mapsto X'$ is a bijection   from~$\cS(\cA)$ onto~$\cS(\cA')$. This bijection is naturally extended to a  bijection between the $\cA$- and $\cA'$-sets, with respect to which the $\cA$-subgroups go to $\cA'$-subgroups. Moreover, if the groups $G$ and $G'$ are abelian, then the last bijection can be extended to the bijection between the $\cA$- and $\cA'$-sections. The images of an $\cA$-set~$X$ and an $\cA$-section~$\fs$ under the corresponding bijections are denoted by $X^\varphi$  and $\fs^\varphi$, respectively. For any such $\fs$, the algebraic isomorphism~$\varphi$ induces an algebraic isomorphism
$$
\varphi_\fs:\cA^{}_{\fs^{}}\to\cA'_{\fs'},\ \pi_{\fs^{}}(X)\mapsto\pi_{\fs'}(X'),
$$
where $\fs'=\fs^\varphi$.

The above bijection between the $\cA$- and $\cA'$-sets is in fact an isomorphism of the corresponding lattices.  Furthermore,  given an $\cA$-set $X$, we have
\qtnl{150114e}
\grp{X^\varphi}=\grp{X}^\varphi\qaq\rad(X^\varphi)=\rad(X)^\varphi.
\eqtn

Let $f:G\to G'$ be a bijection taking the identity of~$G$ to that of $G'$. Then $f$ is a  normalized isomorphism from $\cA$ to $\cA'$ if and only if there exists an algebraic isomorphism~$\varphi:\cA\to\cA'$  such that $X^\varphi=X^f$ for all $X\in\cS(\cA)$. In this case, we say that $f$ \emph{induces} $\varphi$.  The set of all normalized isomorphisms  that induce  a fixed~$\varphi$ is  denoted by  $\iso(\cA,\cA',\varphi)$; we abbreviate
$\iso(\cA,\varphi)=\iso(\cA,\cA,\varphi)$.

\subsection{Circulant S-rings}
An S-ring $\cA$ over a cyclic group $G$ is said to be \emph{circulant}. Any basic set of $\cA$, containing a generator of~$G$, is said to be \emph{highest}. The class of circulant S-rings is closed with respect to tensor product (if the underlying groups of the factors are of coprime orders)  and wreath product.

As was shown in~\cite{Muz1994}, any two algebraically isomorphic circulant S-rings are isomorphic, and are equal if they have the same underlying groups. Therefore, in what follows, we focus on algebraic automorphisms and isomorphisms of a circulant S-ring to itself, rather than general algebraic and combinatorial isomorphisms.  The fact that any two nonisomorphic subgroups of a cyclic group $G$ have different orders, implies that if  $\cA$ is an S-ring over $G$ and $H$ is an $\cA$-subgroup, then $H=H^\varphi=H^f$ for all algebraic automorphisms~$\varphi$ of~$\cA$ and all normalized isomorphisms $f$ of~$\cA$ to itself.

The general theory of circulant S-rings (see, for example, \cite{EvdP2003}) shows that each such ring is obtained by generalized wreath and tensor products from trivial and normal  circulant rings (the normal S-rings are defined below). For the purposes of this paper, there is no need to give this theory in detail. Instead, we define in Subsection~\ref{100924g} a class of totally normal S-rings that is contained in the intersection of the classes of quasinormal and quasidense rings  studied in papers~\cite{EvdP2004a,Wu2024} and~\cite{Evdokimov2016}, respectively. This will allow us to freely use the results of those papers without further specification.

\subsection{Coset closure}\label{250824a} An $S$-ring $\cA$ over a group $G$ is called a \emph{coset} S-ring if every its basic set  is a coset by some $\cA$-subgroup of~$G$. The class of coset S-rings contains all group rings and is closed with respect to tensor and wreath products.

Let $\cA$ be a \emph{quasidense} circulant  S-ring, i.e.,  $\cA$ satisfies the following condition: if~$\fs$ is an $\cA$-section and the S-ring $\cA_\fs$ is trivial, then the group $\fs$ is of prime order. In accordance with paper~\cite{Evdokimov2016}, the \emph{coset closure} $\cA_0$ of $\cA$ is defined to be the intersection of all coset S-rings over $G$, that contain~$\cA$. The S-ring $\cA_0$ is a coset one and $\fS(\cA_0)=\fS(\cA)$, see Theorems~8.5 and~10.7 \emph{ibid}. Denote by $\fS_0(\cA)$ the set  of all $\cA$-sections $\fs$ such that $\fs$ is a subsection of some  $\cA$-section projectively equivalent to a principal $\cA$-section. Then
\qtnl{230824a}
\fS_0(\cA)\subseteq \{\fs\in\fS(\cA):\ (\cA_0)_\fs=\mZ\fs\}=\fS_0(\cA_0),
\eqtn
see formula~(37) and Corollary~10.10 \emph{ibid}.

\subsection{Normal and totally normal S-rings}\label{100924g}
An S-ring $\cA$ over a group $G$ is said to be \emph{normal} if the subgroup $G_{right}$ is normal in~$\aut(\cA)$. The class of normal S-rings contains all group rings and is closed with respect to extensions and tensor products. A complete classification of normal circulant S-rings was obtained in~\cite{EvdP2003}.  In particular, it was proved there  that every normal circulant S-ring is cyclotomic and the following lemma holds true.

\lmml{280724b}
Let $\cA$ be a normal circulant S-ring over a group~$G$ of odd order. Then
\nmrt
\tm{1} $\WL_G(\und{X})=\cA$ for any highest basic set $X \in\cS(\cA)$,
\tm{2} $\iso(\cA)\le\aut(G)$,
\tm{3} $\cA=\cyc(M,G)$ for a unique $M\le\aut(G)$, and also $\aut(\cA)=G\rtimes M$.
\enmrt
\elmm
\prf
The assumption that the group~$G$ is of odd order implies that, in terms of the paper~\cite{EvdP2003}, the S-ring $\cA$ has trivial radical. Thus the first two statements follow from Lemma~6.7 of that paper and the proof of it, whereas  the third statement is a consequence of Theorem~6.1 and Corollary~6.3~\emph{ibid}.
\eprf

An $\cA$-section $\fs$ of a circulant S-ring $\cA$ is said to be \emph{normal} if the S-ring $\cA_\fs$ is normal. Any $\cA$-section projectively equivalent to a normal section is normal. However, in general, a  subsection of a normal section is not necessarily normal.

A circulant S-ring $\cA$ is said to be \emph{totally normal}  if any subsection of a principal $\cA$-section is normal, or, equivalently, any section belonging to~$\fS_0(\cA)$ is normal. Obviously, any circulant coset S-ring is  totally normal, because the restriction of it to any principal section is a group ring (which is always normal). In general, an $\cA$-section $\fs$ of a circulant S-ring $\cA$ is said to be totally normal if the S-ring $\cA_\fs$ is totally normal.

Any totally normal circulant S-ring $\cA$ is quasidense. Indeed,  otherwise, there is an $\cA$-section $\fs$ of composite order and such that the S-ring $\cA_\fs$ is trivial. But then $\fs$ is a subsection of a principal section. Since $\fs$ is not normal, we arrive at a contradiction.

From the above mentioned classification of normal circulant S-rings, it follows that any nontrivial circulant  S-ring over a group of prime order is also totally normal. Note also that the class of totally normal circulant S-rings is closed with respect to tensor products such that the underlying groups of the factors are of coprime orders.

\subsection{Cayley schemes and S-rings}\label{110924a}
Let $G$ be a group. A coherent configuration $\cX=(\Omega,S)$ is called a \emph{Cayley scheme} over~$G$ if
$$
\Omega=G\qaq G_{right}\le\aut(\cX).
$$
Each relation $s \in S$ is the arc set of the Cayley graph  $\cay(G,X)$, where $X=X(s)$  is the neighborhood of the identity of $G$ in the relation~$s$. The subsets $X(s)$, $s\in S$, form a partition $\cS$ of the group~$G$. Moreover, this partition satisfies conditions~(S1), (S2), and (S3). Thus, every Cayley scheme over~$G$ is associated with an S-ring~$\cA=\cA(\cX)$ over~$G$, such that $\cS(\cA)=\cS$.

Conversely, let $\cA$ be an S-ring over $G$. For each $X\in\cS(\cA)$, denote by $s_X$ the arc set of the Cayley graph $\cay(G,X)$. Then the relations $s_X$, $X\in\cS(\cA)$, form a partition~$S$ of the Cartesian square $\Omega^2$, where $\Omega=G$. Moreover, this partition satisfies conditions~(CC1), (CC2), and (CC3). Thus, every S-ring~$\cA$  over~$G$ is associated with a Cayley scheme~$\cX=\cX(\cA)$ over~$G$, such that $S(\cX)=S$.

For any group $G$, the above mappings
\qtnl{140924a}
\cX\mapsto\cA(\cX)\qaq \cA\mapsto\cX(\cA)
\eqtn
define a one-to-one correspondence between the Cayley schemes over $G$ and the S-rings over~$G$. This correspondence preserves the corresponding partial orders on the schemes and S-rings. Moreover, it is functorial in the sense that for any Cayley schemes $\cX$ and $\cX'$, we have
$$
\iso(\cX,\cX')=\iso(\cA,\cA'),
$$
where $\cA=\cA(\cX)$ and $\cA'=\cA(\cX')$. Furthermore, there is a one-to-one correspondence between algebraic isomorphisms, in which an algebraic isomorphism $\varphi:\cX\to\cX'$, $s\mapsto s'$, corresponds to the algebraic isomorphism $\tilde\varphi:\cA\to\cA'$, $X(s)\mapsto X(\varphi(s))$. In what follows, we denote $\tilde\varphi$ just by $\varphi$.

We complete the subsection by remarking that the mappings~\eqref{140924a} respect the tensor and wreath products, namely,
$$
\cA(\cX_1\otimes\cX_2)=\cA(\cX_1)\otimes\cA(\cX_2)\qaq \cA(\cX_1\wr\cX_2)=\cA(\cX_1)\wr\cA(\cX_2).
$$

\section{Multipliers of totally normal circulant S-rings}\label{101024a}
Throughout this section, $\cA$ is a totally normal circulant  S-ring over a (cyclic)  group $G$, $\cS=\cS(\cA)$ and $\fS_0=\fS_0(\cA)$. We denote by $\cA_0$ the coset closure of $\cA$ and put $\cS_0=\cS(\cA_0)$. Finally, we denote by $\iso_\cA(\cA_0)$ the set of all isomorphisms $f\in\iso(\cA_0)$ such that $\cA^f=\cA$, and put
$$
\aut_\cA(\fs)=\aut(\cA_\fs)\cap\aut(\fs)
$$
for each $\fs\in\fS_0$.

\subsection{Outer multipliers}\label{300724a} Suppose that for every section $\fs\in\fS_0$, we are given a coset $\fc_\fs=\sigma_\fs\cdot \aut_\cA(\fs)=\aut_\cA(\fs)\cdot \sigma_\fs$, where $\sigma_\fs\in\aut(\fs)$. Following \cite{Evdokimov2015}, the family
\qtnl{310324a}
\fc=\{\fc_\fs:\ \fs\in\fS_0\}
\eqtn
is called an \emph{outer $\cA$-multiplier} if the following  conditions are satisfied for all
sections $\fs,\ft\in\fS_0$:
\qtnl{230824f}
\fs\succeq \ft\ \Rightarrow\ (\fc_\fs)^\ft=\fc_\ft\quad\qaq\quad
\fs\sim \ft\ \Rightarrow\ (\fc_\fs)^{f_{\fs,\ft}}=\fc_\ft,
\eqtn
where $(\fc_\fs)^\ft=\{f^\ft:\ f\in\fc_\fs\}$ and $f_{\fs,\ft}$ is the projective isomorphism from $\fs$ onto $\ft$. The restriction $\fc_{\fT}$ of the outer multiplier~$\fc$ to a subset $\fT\subseteq\fS_0$ is defined to be the family $\{\fc_\fs:\ \fs\in \fT\}$.

The set of all outer $\cA$-multipliers forms a group $\omult(\cA)$ with respect to  componentwise multiplication.  By \cite[Theorem~3.4]{Evdokimov2015}, there is an isomorphism $\fc\mapsto \varphi_\fc$, $\fc\in \omult(\cA)$, from the group~$\omult(\cA)$ onto the group of all algebraic automorphisms of the S-ring~$\cA$. Here the algebraic automorphism $\varphi=\varphi_\fc$ is defined as follows. Let $X\in\cS(\cA)$ and let $\fs=\fs_X$ the principal section associated with~$X$. Then we put~$X^\varphi$  to be the $\pi_\fs$-preimage of the set $(X_\fs)^{\sigma_\fs}$, see formula~\eqref{240824h}. Thus,
\qtnl{010824a}
(X_\fs)^{\varphi_\fs}=(X_\fs)^{\fc_\fs}=(X_\fs)^{\sigma_\fs}.
\eqtn

\subsection{Inner multipliers}\label{180225a}   Suppose that for every section $\fs\in\fS_0$, we are given an element $\fm_\fs\in\aut(\fs)$.  Following \cite{Evdokimov2015}, the family
\qtnl{310324c}
\fm=\{\fm_\fs:\ \fs\in\fS_0\}
\eqtn
is called an \emph{inner  $\cA$-multiplier} if the conditions~\eqref{230824f} are satisfied for the cosets~$\fc_\fs$ replaced by the singletons $\{\fm_\fs\}$. The restriction $\fm_{\fT}$ of the inner multiplier~$\fm$ to a subset $\fT\subseteq\fS_0$ is defined to be the family $\{\fm_\fs:\ \fs\in \fT\}$.

The set of all inner $\cA$-multipliers forms a subgroup $\mult(\cA)$  of the direct product of the groups $\aut(\fs)$, $\fs\in\fS_0$. For a subgroup $\fM\le \fM(\cA)$ and a set $\fT\subseteq\fS_0$, we put $\fM_\fT=\{\fm_\fT:\ \fm\in\fM\}$ and abbreviate $\fM_\fs:=\fM_{\{\fs\}}$ for $\fs\in\fS_0$. Clearly, $\fM_\fT$ is a group with respect to componentwise product.  For any inner $\cA$-multiplier $\fm$ and   any outer $\cA$-multiplier $\fc$, we write $\fm_\fT\le\fc_\fT$ if $\fm_\fs\in\ \fc_\fs$ for all $\fs\in\fT$. When $\fT=\fS_0$, we abbreviate $\fm\le\fc$.

\lmml{310324f}
For each  inner $\cA$-multiplier $\fm\in\mult(\cA)$, there exists a normalized isomorphism $f_\fm \in\iso_\cA(\cA_0)$ such that for all  $X\in\cS$, we have
\qtnl{160524a}
{X_\fs}^{(f_\fm)^\fs}={X_\fs}^{\fm_\fs},
\eqtn
where $\fs=\fs_X$.
\elmm
\prf
By \cite[Lemma~5.2]{Evdokimov2015}, the family of the singletons $\{\fm_\fs\}$, $\fs\in\fS_0$, is an outer $\cA_0$-multiplier. By \cite[Theorem~3.4]{Evdokimov2015},  there exists a uniquely determined algebraic automorphism $\varphi$ of the S-ring~$\cA_0$ such that
\qtnl{010423a}
{X_\fs}^{\varphi_\fs}={X_\fs}^{\fm_\fs}, \qquad  X\in\cS,
\eqtn
where $\fs=\fs_X$. By \cite[Theorem 9.1]{Evdokimov2016}, every algebraic automorphism of $\cA_0$ is induced by an isomorphism. Consequently, there exists a normalized isomorphism $f\in\iso(\cA_0,\varphi)$. By virtue of~\eqref{010423a}, we have
$$
{X_\fs}^{f^\fs}={X_\fs}^{\varphi_\fs}={X_\fs}^{\fm_\fs},
$$
which proves formula~\eqref{160524a} for $f_\fm=f$.

To prove that  $f\in\iso_\cA(\cA_0)$, we need to verify that $X^f\in\cS$ for all $X\in\cS$.  The principal section $\fs=\fs_X$ belongs to the set $\fS_0$. By formula~\eqref{240824h},
$$
X^f=(\pi^{-1}(X_\fs))^f=\pi^{-1}({X_\fs}^{f^\fs})=\pi^{-1}({X_\fs}^{\fm_\fs}),
$$
where $\pi=\pi_\fs$. Put $Y=X^{\fm'}$, where  $\fm'\in\aut(G)$ is such that $(\fm')^{\fs}=\fm_\fs$. By Theorem~\ref{261009b}(1), we have $Y\in\cS$. Furthermore, $\grp{X}=\grp{Y}$ and $\rad(X)=\rad(Y)$. Therefore, $\fs=\fs_X=\fs_Y$ and $Y_\fs={X_\fs}^{\fm_\fs}$. Thus,
$$
X^f=\pi^{-1}({X_\fs}^{\fm_\fs})=\pi^{-1}({Y_\fs})=Y\in\cS,
$$
and $\cA^f=\cA$, as required.
\eprf

Each normalized isomorphism $f\in\iso(\cA)$ leaves every $\cA$-subgroup fixed (because the S-ring $\cA$ is circulant). This yields the restriction isomorphism $f^\fs\in\iso(\cA_\fs)$ for every $\cA$-section~$\fs$. Note that the collection $\fm_f=\{f^\fs:\ \fs\in\fS_0\}$ is not necessarily an  inner $\cA$-multiplier (because $f^\fs$ does not always belong to $\aut(\fs)$). Put
$$
\mult_{aut}(\cA)=\{\fm_f\in\mult(\cA):\ f\in\aut(\cA)\}.
$$
Assume that the group $G$ is of odd order. Let $\fs\in\fS_0$. Since the S-ring $\cA$ is totally normal, the S-ring $\cA_\fs$ is normal  and by Lemma~\ref{280724b}(2), the isomorphism $f^\fs$ belongs to $\aut(\fs)$. Thus, $\fm_f\in\mult(\cA)$ and $\mult_{aut}(\cA)$ is a subgroup of $\mult(\cA)$.

\subsection{Algebraic fusions}
Let $\cA$ be a coset S-ring, i.e., $\cA=\cA_0$, and  let $\fM$ be a group of inner $\cA$-multipliers.  This group acts on the basic sets of $\cA$, namely,  if $X\in\cS$ and $\fm\in \fM$, then
\qtnl{040824j1}
X^\fm:=\pi_\fs^{-1}({X^{}_\fs}^{\fm_\fs}),
\eqtn
where $\fs=\fs_X$. We define the \emph{algebraic fusion} of $\cA$ with respect to the group $\fM$ as an S-ring $\cA^\fM$ over  $G$ with  the basic sets
\qtnl{040824j2}
X^\fM=\bigcup_{\fm\in\fM}X^\fm, \qquad X\in \cS.
\eqtn
In the sense of \cite[Subsection~2.3.3]{CP2019}, the Cayley scheme $\cX(\cA^\fM)$ is the algebraic fusion of the Cayley scheme $\cX(\cA)$ with respect to the group of the algebraic automorphisms $\varphi_{\fc_\fm}$ corresponding to the outer $\cA$-multiplier $\fc_\fm=\{\{\fm_\fs\}:\ \fs\in\fS_0\}$, $\fm\in\fM$. Since the S-ring $\cA$ is a coset one, formula~\eqref{230824a} implies that
\qtnl{040824g}
(\cA^\fM)_\fs=(\cA_\fs)^{\fM_\fs}=(\mZ\fs)^{\fM_\fs}=\cyc(\fM_\fs,\fs).
\eqtn
 for all $\fs\in\fS_0$.

Obviously, $\cA^\fM\le\cA$. Therefore each $\cA^\fM$-subgroup is an $\cA$-subgroup. Conversely, if $H$ is an $\cA$-subgroup, then by Lemma~\ref{310324f}, we have $H^\fm=H^{f_\fm}=H$ for all $\fm\in\fM$, and hence $H$ is  an $\cA^\fM$-subgroup. Thus  the S-rings $\cA$ and $\cA^\fM$ have the same $\cA$-subgroups and
\qtnl{040824q}
\fS(\cA)=\fS(\cA^\fM).
\eqtn
Moreover, it is not hard to see that the relations $\preceq$ and $\sim$  on the $\cA$- and $\cA^\fM$-sections are the same.

\lmml{120724a}
In the above notation,  assume that $G$ is of odd order and the S-ring $\cyc(\fM_\fs,\fs)$ is  normal for all $\fs\in\fS_0(\cA)$. Then $\fM\le\fM_{aut}(\cA^\fM)$ and
\qtnl{120724b}
\fS_0(\cA^\fM)=\fS^{}_0(\cA)\qaq(\cA^\fM)_0=\cA.
\eqtn
\elmm
\prf
To prove the first inclusion, we need to verify that for each $\fm\in\fM$, there exists $f\in\aut(\cA^\fM)$ such that $\fm=\fm_f$.  To this end, denote by $f$ the isomorphism $f_\fm\in\iso(\cA)$ from Lemma~\ref{310324f} (note that in our case, $\cA=\cA_0$). Then $\fm=\fm_f$. Moreover, for  all $X\in\cS$, we have
$$
(X^\fM)^f=\Bigl(\bigcup_{\fn\in\fM}X^\fn\Bigr)^f
=\Bigl(\bigcup_{\fn\in\fM}\pi_\fs^{-1}({X_\fs}^{\fn_\fs})\Bigr)^f=
\bigcup_{\fn\in\fM}\pi_\fs^{-1}({X_\fs}^{\fn_\fs\fm_\fs})
$$
$$
=\bigcup_{\fn\in\fM}\pi_\fs^{-1}({X_\fs}^{(\fn\fm)_\fs})=
\bigcup_{\fn'\in\fM}\pi_\fs^{-1}({X_\fs}^{(\fn')_\fs})=\bigcup_{\fn'\in\fM}X^{\fn'}=X^\fM,
$$
where $\fs=\fs_X$, see formulas~\eqref{040824j1} and \eqref{040824j2}. Thus, $f\in\aut(\cA^\fM)$, as required.

First, we prove that $\fS_0(\cA)\subseteq \fS_0(\cA^\fM)$.  Since the relations $\preceq$ and $\sim$  on the $\cA$- and $\cA^\fM$-sections are the same, it suffices to verify that every  principal $\cA$-section $\fs$ is a  principal  $\cA^\fM$-section. Suppose that $\fs=\fs_{X}$ for some $X\in\cS$. Then obviously $\grp{X}=\grp{X^\fM}$ and $\rad(X)\le\rad(X^\fM)$. Suppose on the contrary that $\rad(X)\neq \rad(X^\fM)$. Then the radical of the set  $(X^\fM)_\fs\le \fs$ is nontrivial. Since the group $\fs$ is of odd order and $(X^\fM)_\fs$ is a highest basic set of the S-ring~$(\cA^\fM)_\fs=\cyc(\fM_\fs,\fs)$, this S-ring is not normal by \cite[Theorem~6.1]{EvdP2003}, a contradiction.

From formula~\eqref{230824a}, it follows that $\fS_0(\cA^\fM)\subseteq\fS_0((\cA^\fM)_0)$. On the other hand, $\cA$ is a coset S-ring containing $\cA^\fM$. Therefore, $ (\cA^\fM)_0\le \cA$. It follows that $\fS_0((\cA^\fM)_0)\subseteq\fS_0(\cA)$. Consequently, $\fS_0(\cA^\fM)\subseteq \fS_0(\cA)$. Together with the inclusion in the above paragraph, this proves the first equality in~\eqref{120724b}. Furthermore, the above argument shows that $\fS_0((\cA^\fM)_0)=\fS_0(\cA)$. Consequently, the coset S-rings $(\cA^\fM)_0$ and $\cA$ have the same principal sections. Thus they are equal which proves the second equality in~\eqref{120724b}.
\eprf

\section{Local system of multipliers }\label{101024b}

Let $\cA$ be a circulant S-ring, $\fc$ an outer $\cA$-multiplier, and $k\ge 2$ an integer. Denote by $\fS_k=\fS_k(\cA)$ the set of all subsets consisting of no more than~$k$ principal $\cA$-sections. Assume that for each $\fS\in \fS_k$, we are given  an inner $\cA$-multiplier~$\fm(\fS)$. A family
\qtnl{160524s}
\{\fm(\fS)\in\fM(\cA):\ \fS\in \fS_k\}
\eqtn
is called a \emph{$k$-local system of inner $\cA$-multipliers} with respect to  $\fc$ if the following two conditions are satisfied  for all subsets $\fS,\fS'\in\fS_k$ and some inner $\cA$-multipliers $\fm(\fS,\fS')\in\mult_{aut}(\cA)$:
\qtnl{080423a}
\fm(\fS)_\fS\le \fc_\fS\qaq \fm(\fS)_{\fS\cap \fS'}=\fm(\fS')_{\fS\cap\fS'}\cdot \fm(\fS,\fS')_{\fS\cap\fS'}.
\eqtn
The main result of this section is  the following theorem.

\thrml{310324ab}
Let $\cA$ be a totally normal circulant S-ring and $k\ge 2$ an integer. Assume that
\nmrt
\tm{1} the scheme $\cX(\cA_0)$ is binary,
\tm{2} there is a $k$-local system of inner $\cA$-multipliers with respect to  an outer multiplier $\fc\in\fC(\cA)$.
\enmrt
Then for each $2\le m< \sqrt{k}$, the scheme   $\cX(\cA)$  is $\WL_m$-equivalent to itself with respect to  the algebraic isomorphism $\varphi_\fc$.
\ethrm

We prove Theorem~\ref{310324ab} at the end of the section. In what follows, $G$ is a cyclic group underlying the S-ring $\cA$,  $\cX=\cX(\cA)$ is the Cayley scheme associated with $\cA$, and $\varphi=\varphi_\fc$. Fix an integer $\ell$ such that
\qtnl{240824f}
2\le \ell\le \lfloor\sqrt{k}\rfloor,
\eqtn
and fix a $k$-local system~\eqref{160524s} of inner $\cA$-multipliers with respect to  the outer multiplier~$\fc$.

Denote by~$\cR_\ell=\cR_\ell(\cX)$ the partition of the Cartesian power $G^\ell$ into the classes of the equivalence relation $\sim_\ell$ defined as follows: $x \sim_\ell y$ if and only if $R_\ell(x)=R_\ell(y)$, where  $R_\ell(\cdot)=R_{\cX,\ell}(\cdot)$ is the $\ell\times \ell$-array defined in formula~\eqref{170423a}. In particular, each class $\sX\in \cR_\ell$ is a maximal subset of $G^\ell$ such that the array $R_\ell(x)$  does not depend on the choice of the $\ell$-tuple $x\in \sX$. In accordance with \cite[Subsection~3.1]{Chen2023a}, we have
\qtnl{160524a1}
\WL_\ell(\cX)\ge \cR_\ell(\cX),
\eqtn
i.e., each class of the partition $\cR_\ell(\cX)$ is a union of some classes of the partition~$\WL_\ell(\cX)$. Furthermore, the algebraic isomorphism $\varphi$ induces a permutation of $\cR_\ell(\cX)$, that takes $\sX$ to $\sX^\varphi=\{y\in G^\ell:\ R(y)=R(x)^\varphi\}$, where $R(x)^\varphi=(r(x_i,x_j)^\varphi)$ and  $x\in \sX$. In particular,
\qtnl{250525a}
\cR_\ell(\cX)^{\hat\varphi}=\cR_\ell(\cX).
\eqtn

Let $\sX$ be a class of the partition $\cR_\ell$, and let $1\le i,j\le \ell$. The basis relation $s=r(x_i,x_j)$ of the  scheme~$\cX$ does not depend on the $\ell$-tuple $x\in \sX$. Denote by $X_{ij}$ the basic set of the S-ring~$\cA$, corresponding to $s$. Clearly,   $x_j{x_i}^{-1}\in X_{ij}$ for all $x\in\sX$. The set
\qtnl{120225a}
\fS(\sX)=\{\fs_{X_{ij}}:\  1\le i,j\le \ell\}
\eqtn
consists of at most $\ell^2\le k$ principal $\cA$-sections  $\fs_{X_{ij}}$ and hence $\fS(\sX)\in\fS_k$.  Denote by $\fm(\sX)$  the inner $\cA$-multiplier~$\fm(\fS)$  belonging to the family~\eqref{160524s} for $\fS=\fS(\sX)$, and by $f_\sX$  the isomorphism $f_{\fm(\sX)}\in\iso_\cA(\cA_0)$ in Lemma~\ref{310324f} for $\fm=\fm(\sX)$.

\lmml{160524f}
In the above notation, the mapping
$$
f:G^\ell\to G^\ell,\   x\mapsto x^{f_\sX},
$$
is a bijection, where $\sX\in \cR_\ell(\cX)$ is the class containing the $\ell$-tuple $x$ and $x^{f_\sX}$ denotes the $\ell$-tuple $(x_1^{f_\sX},\ldots,x_\ell^{f_\sX})$. Moreover,
\qtnl{170524f}
\cR_\ell(\cX)^f=\cR_\ell(\cX)\qaq\WL_\ell(\cX)^f=\WL_\ell(\cX).
\eqtn
\elmm
\prf
By \eqref{250525a}, to prove that $f$ is a bijection and the first equality holds true, it suffices to verify that for any$\sX\in\cR_\ell(\cX)$, we have  $\sX^f=\sX^\varphi$. To do this, we need to show that $R_\ell(x^f)=R_\ell(x)^\varphi$, or, equivalently, that  $r({x_i}^f,{x_j}^f)=r(x_i,x_j)^\varphi$ for all   $1\le i,j\le \ell$.  However, since $f_\sX$ is a normalized isomorphism of the S-ring~$\cA$, we have
$$
x_j{x_i}^{-1}\in X_{ij}\quad\Leftrightarrow\quad {x_j}^{f_\sX}({x_i}^{f_\sX})^{-1}\in (X_{ij})^{f_\sX}.
$$
By the first condition in~\eqref{080423a}, we have $(X_{ij})^{f_\sX}={X_{ij}}^\varphi$, and we are done.

To prove the second equality in~\eqref{170524f}, we need to verify that $\sX^f\in\WL_\ell(\cX)$ for each $\sX\in\WL_\ell(\cX)$. Note that by virtue  of~\eqref{160524a1}, there is a class $\sX'\in\cR_\ell(\cX)$ containing the class~$\sX$. Therefore, $\sX^f=\sX^{f_{\sX'}}$. Since $f_{\sX'}\in\iso(\cX)$, we conclude that $\sX^{f_{\sX'}}\in\WL_\ell(\cX)$. Thus, $\sX^f\in\WL_\ell(\cX)$,  as required.
\eprf

The following lemma is a key point in the proof of Theorem~{\rm \ref{310324ab}}. Below, for an $\ell$-tuple $x$, we put  $x\cdot \alpha= (x_1,\ldots,x_\ell,\alpha)$ for all $\alpha\in G$. For any set $\sX\subseteq G^{\ell+1}$, we put $\res_x \sX=\{\alpha\in G:\ x\cdot \alpha\in \sX\}$.

\lmml{170524z}
Assume that $m<\lfloor \sqrt{k}\rfloor$. Then given a tuple $x\in G^m$, there is a permutation  $\theta\in\sym(G)$  such that  for all $\alpha\in G$
\qtnl{180524a}
R(x^f\cdot \alpha^\theta)=R(x\cdot \alpha)^\varphi,
\eqtn
where $R=R_{m+1}$ and $f$ is the bijection from Lemma~{\rm \ref{160524f}} for $\ell=m$.
\elmm
\prf
Denote by $\hat f:G^{m+1}\to G^{m+1}$ the bijection $f$ defined in Lemma~\ref{160524f} for  $\ell=m+1$; note that by the assumption on $m$, we have $\ell\le\lfloor \sqrt{k}\rfloor$. It suffices to prove that if $\hat \sX\in \cR_{m+1}(\cX)$ 
and $x\in\pr_m \hat \sX$, then there exists a bijection $\theta_{\hat \sX}:\res_x \hat \sX\to \res_{x^f} \hat \sX^{\hat f}$ such that
\qtnl{190724a}
R(x^f\cdot \alpha^{\theta_{\hat \sX}})=R(x\cdot \alpha)^\varphi \quad\text{for all }\alpha\in \res_x \hat \sX.
\eqtn
Indeed,  $G$ is  a disjoint union of  
those sets $\res_x \hat \sX$ for which  $\hat \sX\in \cR_{m+1}(\cX)$ and $x\in\pr_m \hat \sX$. Moreover, assuming the existence of the bijections~$\theta_{\hat \sX}$, we see that
$$
\res_{x^f} \hat \sX^{\hat f}\ne\varnothing \quad \Leftrightarrow\quad \res_x\hat \sX\ne\varnothing.
$$
Thus the required bijection $\theta:G\to G$ can be defined as follows: $\theta(\alpha)=\alpha^{\theta_{\hat \sX}}$, where $\hat\sX$ is the class containing the tuple $x\cdot\alpha$.

To prove formula~\eqref{190724a}, fix $\alpha\in \res_x \hat \sX$ and put $\sX=\pr_m \hat \sX$. Let $\hat\fS=\fS(\hat \sX)$ and $\fS=\fS(\sX)$ be the sets of principal $\cA$-sections, defined by formula~\eqref{120225a}.   Then $\fS\subseteq \hat\fS$, and $\fS,\hat\fS\in \fS_k$ by the choice of $m$. Since $x\in \sX$ and $x\cdot \alpha\in\hat \sX$, this implies that
\qtnl{120225t}
x^f=x^{f_\sX}\qaq (x\cdot \alpha)^{\hat f}=(x\cdot \alpha)^{\hat f_{\hat \sX}}=x^{\hat f_{\hat \sX}}\cdot \alpha^{\hat f_{\hat \sX}}.
\eqtn
According to our definition,   $f_\sX=f_{\fm(\sX)}$ and $\hat f_{\hat \sX}=f_{\fm(\hat\sX)}$ are the isomorphisms of the S-ring~$\cA$, inducing the inner $\cA$-multipliers $\fm(\fS)$ and $\fm(\hat\fS)$, respectively. Now the left-hand side of formula~\eqref{080423a} implies that
$$
\hat f^{\fs_{ij}}=\fm(\fS')_{\fs_{ij}}\in \fc_{\fs_{ij}}
$$
for all principal $\cA$-sections $\fs_{ij}=\fs_{X_{ij}}$ belonging to the subset~$\hat\fS$. This implies that $(X_{ij})^{\hat f}=\varphi(X_{ij})$. By the right-hand side of formula~\eqref{120225t}, this yields
\qtnl{190424e}
R(x^{\hat f}\cdot \alpha^{\hat f})=R(x\cdot \alpha)^\varphi\quad \text{for all } \alpha\in \res_x \hat \sX.
\eqtn
At this point, we need an auxiliary statement concerning the array $R_0(\cdot)=R_{\cX_0,m}(\cdot)$ defined by formula~\eqref{170423a}, where $\cX_0=\cX(\cA_0)$ (recall that $\cA_0$ is the coset closure of~$\cA$).

\clml{bb}
Let $h=h(\fS,\hat\fS)$ be the automorphism of the S-ring $\cA$, corresponding to the inner $\cA$-multiplier~$\fm(\fS,\hat\fS)\in\fM_{aut}(\cA)$ on the right-hand side of formula~\eqref{080423a}. Then
$$
R_0(x^{\hat f})=R_0(x^{fh}).
$$
\eclm
\prf
We need to verify that for all $1\le i,j\le m$, the relations  $r_0((x^{\hat f})_i,(x^{\hat f})_j)$ and $r_0((x^{fh})_i,(x^{fh})_j)$ are equal, where $r_0=r_{\cX_0}$. Recall that $\hat f_{\hat \sX}\in\iso_\cA(\cA_0)$ and hence  $\hat f_{\hat \sX}\in\iso(\cX_0)$. Therefore,
\qtnl{020824a}
r_0((x^{\hat f})_i,(x^{\hat f})_j)=r_0({x_i}^{\hat f_{\hat\sX}},{x_j}^{\hat f_{\hat\sX}})=
r_0(x_i,x_j)^{\hat f_{\hat\sX}}=r_0(x_i,x_j)^{\hat f}.
\eqtn
 Similarly, since $f_\sX$ and $h_\sX$ are also isomorphisms of the scheme~$\cX_0$, we have
\qtnl{020824a1}
r_0((x^{fh})_i,(x^{fh})_j)=r_0(x_i,x_j)^{fh}.
\eqtn
It remains to verify that the relations on the right-hand sides of formulas~\eqref{020824a} and~\eqref{020824a1} coincide. Denote by $X_0$  the basic set  of the S-ring $\cA_0$ that corresponds to the basis relation $r_0(x_i,x_j)$. Then $X_0$  contains the element $x_0=x^{}_jx_i^{-1}$. The element $\bar x_0=(x_0)_\fs$ belongs to the set $X_{ij}$ that generates the cyclic group $\fs=\fs_{X_{ij}}$. Since the S-ring $\cA$ is totally normal, the S-ring $\cA_\fs$ is normal and hence every element of~$X_{ij}$ is a generator of~$\fs$; in particular, $\fs=\grp{x_0}$. It follows that $\fs=\fs_{X_0}$. Consequently,
$$
{\bar x_0\,}^{(\hat f_{\hat \sX})^\fs}={\bar x_0\,}^{ \hat\fm_\fs}={\bar x_0\,}^{ \fm_\fs\fn_\fs}={\bar x_0\,^{(f_\sX)^\fs (h_\sX)^\fs}}={\bar x_0\,}^{(f_\sX h_\sX)^\fs},
$$
where $\hat\fm=\fm(\hat\fS)$, $\fm=\fm(\fS)$, $\fn=\fm(\fS,\hat\fS)$. and the equality $\hat\fm_\fs=m_\fs n_\fs$ follows from the second equality in~\eqref{080423a}. Therefore,
$$
{X_0}^{\hat f}=(\pi^{-1}({\bar x_0}))^{\hat f}=\pi^{-1}({\bar x_0\,}^{(\hat f_{\hat \sX})^\fs})=\pi^{-1}({\bar x_0\,}^{(f_\sX h_\sX)^\fs})=(\pi^{-1}({\bar x_0}))^{fh}={X_0}^{fh},
$$
where $\pi=\pi_\fs$. Thus, ${X_0}^{\hat f}={X_0}^{fh}$, as required.
\eprf

By theorem condition, the scheme $\cX_0$ is binary. Besides, $R_0(x^{\hat f})=R_0(x^{fh})$ by Claim~\ref{bb}, where $h\in \aut(\cA)$. Thus there exists  $h_0\in\aut(\cX_0)\le \aut(\cX)$ such that $x^{\hat f h_0}=x^{fh}$. Then
\qtnl{200724a}
(x^{\hat f}\cdot \alpha^{\hat f})^{h_0}=x^{\hat f h_0}\cdot \alpha^{\hat fh_0}=x^{fh}\cdot \alpha^{\hat fh_0h^{-1}h}=(x^f\cdot \alpha^{\hat fh_0h^{-1}})^h.
\eqtn
Note that $x^{\hat f}\cdot \alpha^{\hat f}\in {\hat \sX}^{\wh f}$ by the choice of~$\alpha$. Since $h_0$ and $h^{-1}$ are  automorphisms of the scheme $\cX$, this implies that  $x^f\cdot \alpha^{\theta_{\hat \sX}}\in {\hat \sX}^{\wh f}$, where $\theta_{\hat \sX}=\hat fh_0h^{-1}$  is  the  composition of isomorphisms of~$\cX$. In particular,
$$
\alpha^{\theta_{\hat \sX}}\in \res_{x^f}\hat \sX^{\hat f}.
$$

To complete the proof, we recall that the point $\alpha$ was  chosen in $\res_x\hat \sX$ arbitrarily. Moreover,  none of the isomorphisms $\hat f$, $h_0$, and~$h$ depends on $\alpha$, and the mapping $\theta_{\hat \sX}:\res_x\hat \sX\to \res_{x^f} \hat \sX^{\hat f}$ is a bijection. Finally, by formulas~\eqref{190424e} and~\eqref{200724a} and because of $h_0,h\in\aut(\cX)$,  we have
$$
R(x\cdot \alpha)^\varphi =R(x^{\hat f}\cdot \alpha^{\hat f})=R(x^{\hat f}\cdot \alpha^{\hat f})^{h_0h^{-1}}=R(x^f\cdot \alpha^{\hat fh_0h^{-1}})=R(x^f\cdot \alpha^{\theta_{\hat \fX}})
$$
which proves formula~\eqref{190724a}.
\eprf

\begin{proof}[Poof of Theorem~{\rm \ref{310324ab}}] By Lemma~\ref{040924b}, to prove that the scheme $\cX$ is $\WL_m$-equiva\-lent to itself with respect to  the algebraic isomorphism~$\varphi$, it suffices to verify that the bijection~$f$ defined in Lemma~\ref{160524f}  for $\ell=m$ satisfies the condition of Lemma~\ref{040924b}. In other words, we need to verify that if $x\in G^m$ and $x'=x^f$, then   Duplicator has a winning strategy in the pebble game  $\sC_{m+1}(\varphi)$ on the scheme~$\cX$ with initial configuration $(x,x')$.  Let Spoiler choose a set~$B\subseteq G$. Duplicator responds with the set $A=B^\theta$, where $\theta$ is the bijection defined in Lemma~\ref{170524z}. Spoiler places one of the pebbles  on a point $\alpha\in A$. Duplicator places the corresponding pebble to the point $b=\alpha^\theta$ and wins, because $R(x\cdot \alpha)^\varphi=R(x'\cdot b)$ by formula~\eqref{180524a}.
\eprf

\section{Two ingredients for the main construction}\label{190524h}
In the sequel, $a=2b$ is a positive  integer and  $V$ is a set of cardinality $a$. We fix a simple  graph $\cG=(V,E)$ belonging to the class $\fB_a$ defined in Section~\ref{120225g}; in particular, $|E|=3a/2$. It is assumed that $V=\{v_1,\ldots,v_a\}$. The vertices  of the parts $V_1$ and $V_2$ of the (bipartite) graph~$\cG$ are denoted by $u_1:=v_1,\ldots,u_b:=v_b$ and $w_1:=v_{b+1},\ldots,w_b:=v_a$, respectively.

The graph $\cG$ being cubic and bipartite admits a $1$-factorization~$F_\cG$, i.e., the set of three pairwise disjoint matchings the union of  which is equal to~$E$. It is convenient to index the matchings  of $F_\cG$ by the nonidentity elements of  the Klein group~$K=C_2\times C_2$; we set $F_\cG=\{E_k:\ k\in K^\#\}$, where $K^\#=K\setminus \{1\}$ and $E_k$ is the matching indexed by~$k$.

\subsection{A coherent configuration associated with $\cG$}\label{121123u}
With each  vertex $v\in V$ of the graph~$\cG$, we associate a copy $\Omega_v$ of the Klein group~$K$ and, for a subset $V'\subseteq V$, we denote by~$\Omega_{V'}$ the union of $\Omega_v$, $v\in V'$. When $V'=\{v,v'\}$, we abbreviate $\Omega_{v,v'}:=\Omega_{\{v,v'\}}$.

For any $v\in V$ and any $k\in K$, we denote by $\bar k_v$ the permutation $x\mapsto xk$, $x\in K$, on the copy $\Omega_v$ of the group~$K$. The set $\bar K_v=\{\bar k_v:\ k\in K\}$ is a permutation group on~$\Omega_v$, that is isomorphic to~$K$. The direct product
$$
K_V=\prod_{v\in V}\bar K_v,
$$
is a permutation group on the set $\Omega=\Omega_V$. It was proved in~\cite[Section~5]{EvdP1999c}  that there exists  a coherent configuration $\cX=\cX(\cG,F_\cG) =(\Omega,S)$ that satisfies the following conditions:
\nmrt
\tm{K1} $F(\cX)=\{\Omega_v:\ v\in V\}$, and also $\aut(\cX)^{\Omega_v}=\bar K_v$ for all $v$,
\tm{K2}  the group of all $f\in\iso(\cX)$ such that ${\Omega_v}^f=\Omega_v$ and $f^{\Omega_v} \in \bar K_v$ for all  $v\in V$, coincides with $K_V$,
\tm{K3} $\aut(\cX)^{\Omega_{v,v'}}=\bar K_{v^{}}\times \bar K_{v'}$ for all distinct $v,v'\in V$ nonadjacent in $\cG$,
\tm{K4}  if  $\{v,v'\}\in E_k$ for some $k\in K^\#$, then
$$
\aut(\cX)^{\Omega_{v,v'}}=\{(x,x')\in \bar K_{v^{}}\times\bar K_{v'}:\ \bar\rho_{v^{}}(x)=\bar\rho_{v'}(x')\},
$$
where $\bar\rho_{v^{}}:\bar K_{v^{}}\to \sym(2)$ and  $\bar\rho_{v'}:\bar K_{v'}\to \sym(2)$ are the natural epimorphisms with kernels $\grp{\bar k_{v^{}}}$ and $\grp{\bar k_{v'}}$, respectively,
\tm{K5} if  $\{u,w\}\in E$ and  $\psi\in\sym(S)$  are such that $\psi(s)\ne s$ if and only if $s\in S_{\Omega_u,\Omega_w}\,\cup\, S_{\Omega_w,\Omega_u}$, then $\psi$ is an algebraic automorphism of~$\cX$.
\enmrt
The details concerning the construction and other properties of the coherent configuration $\cX$ can also be found in~\cite[Subsections~4.1.3 and~4.2.1]{CP2019}, where $\cX$ is defined as a special case of a cubic Klein configuration.

To state the main fact about the coherent configuration $\cX$, that we need later, we fix an edge $e_0\in E$ and set $\psi=\psi_{e_0}$ to be the algebraic isomorphism in condition~(K5) for $e=e_0$. Next, for an integer $m\ge 1$, we put
$$
F_m=F_m(\cX)=\{\Omega_{V'}:\ |V'|\le m\}.
$$
A family $\{f_\Delta\in K_V: \Delta \in F_m\}$ is called an \emph{$m$-local system of isomorphisms for $\cX$} with respect to $\psi$ if for  all $\Delta,\Delta'\in F_m$, we have
$$
(f_\Delta)^\Delta\in\iso(\cX_\Delta,\psi_\Delta)\qaq (f_\Delta {f_{\Delta'}}^{-1})^{\Delta\cap\Delta'}\in\aut(\cX)^{\Delta\cap\Delta'}.
$$

\lmml{050524t}
In the above notation, assume that each separator\footnote{A separator of the graph $\cG$ is a subset $V'\subset V$ satisfying the following condition: each connected component of the graph obtained from $\cG$ by removing all the vertices of $V'$ contains at most half vertices of $\cG$.} of the graph $\cG$ has  cardinality greater than $3m$. Then   $\iso(\cX,\psi)=\varnothing$ and there is  an $m$-local system of isomorphisms for $\cX$ with respect to  $\psi$.
\elmm
\prf
Follows from Lemma~4.2.6 in~\cite{CP2019} and its proof.
\eprf

\subsection{A coset S-ring associated with the graph $\cG$}\label{130524f}
Let $a$ and $b$ be the integers as above, and let $n$ be as in Theorem~\ref{090524d}. We define an  S-ring over the group $G=C_n$ as the tensor product of two iterated wreath products of the group rings of cyclic groups,
\qtnl{121123i}
\cA=(\mZ C_{n_1}\wr \mZ C_{n_2}\wr\cdots \wr \mZ C_{n_b})\otimes (\mZ C_{n_{b+1}}\wr \mZ C_{n_{b+2}}\wr\cdots \wr \mZ C_{n_{2b}}),
\eqtn
where the numbers $n_i$ are as in Theorem~\ref{090524d}.
For any  vertex $v\in V$ of the graph~$\cG$, we denote by $G_v$ the (cyclic) subgroup of $G$ such that
$$
|G_v|=\css
n_1\cdot\ldots \cdot n_i &\text{if $v=u_i$ with $1\le i\le b$,}\\
n_{b+1}\cdot\ldots \cdot n_{b+j} &\text{if $v=w_j$ with $1\le j\le b$.}\\
\ecss
$$
 Clearly, $G=G_{u_b}\times G_{w_b}$ and $\cA$ is the tensor product of  the iterated wreath products $\cA_{G_{u_b}}$ and $\cA_{G_{w_b}}$ of the group rings $\mZ C_{n_i}$. The  statement below immediately follows from remarks in~Subsection~\ref{250824a}) and Theorem~\ref{250724a}.

\lmml{250824b}
The S-ring $\cA$ is a coset one  and the (circulant) scheme $\cX(\cA)$ is binary.
\elmm

We enlarge the sets  $V_1$ and $V_2$ by adding to them new elements $u_0$ and $w_0$, respectively, and set $G_{u_0}=G_{w_0}=1$. Then $G_v$ is an $\cA$-subgroup for each $v$, including these new elements. Moreover,
\qtnl{080724a}
G_{u_0}<G_{u_1}<\ldots<G_{u_b}\qaq G_{w_0}<G_{w_1}<\ldots<G_{w_b}.
\eqtn
We sometimes identify the subgroups $G_{u_i}$ and $G_{u_i}\times G_{w_0}$, and $G_{w_j}$ with $G_{u_0}\times G_{w_j}$. Below, for any vertex $v$ of the graph $\cG$, we define the section
$$
\fs_v=\css
G_{u_i}/G_{u_{i-1}} &\text{if $v=u_i$ with $1\le i\le b$,}\\
G_{w_j}/G_{w_{j-1}} &\text{if $v=w_j$ with $1\le j\le b$.}\\
\ecss
$$

\lmml{101123w}
\phantom{a}
\nmrt
\tm{1} Each $\cA$-subgroup is of the form $U\times W$, where $G_{u_{i-1}}\le U\le G_{u_i}$ and $ G_{w_{j-1}}\le W\le G_{w_j}$ for some $1\le i,j\le b$,
\tm{2} $\fs\in\fS_0(\cA)$ if and only if $\fs\preceq \fs_{u,w}$ for some $u\in V_1$ and $w\in V_2$, where $\fs_{u,w}=\fs_u\times\fs_w$.
\enmrt
\elmm
\prf
Since the group $G=G_{u_b}\times G_{w_b}$ is cyclic, any subgroup $H\le G$ is the direct product of two subgroups $U:=H\cap G_{u_b}$ and $W:=H\cap G_{w_b}$. Assume that $H$ is an $\cA$-subgroup. Since obviously $G_{u_b}$ and $G_{w_b}$ are also $\cA$-subgroups, so is $U$ and hence it is an  $\cA_{G_{u_b}}$-subgroup. By the first part of formula~\eqref{080724a}, there is the smallest $i\in\{1,\ldots,b\}$ such that $U\le G_{u_i}$. Without loss of generality, we assume that $U\ne 1$. Then $U$ is an $(\cA_{i-1}\wr \mZ C_{n_i})$-subgroup, where $\cA_{i-1}$ is the restriction of~$\cA$ to $G_{u_{i-1}}$. The minimality of~$i$ implies that  $U\ge  G_{u_{i-1}}$ (see Subsection~\ref{150225b}). Thus, $G_{u_{i-1}}\le U\le G_{u_i}$. Similarly, one can prove that $ G_{w_{j-1}}\le W\le G_{w_j}$ for some $1\le j\le b$. Statement~(1) is proved.

To prove statement (2), we note that the necessity is a consequence of the fact that $\fs_{u,w}$ is a principal section. Conversely, let $\fs$ be an $\cA$-section. By statement~(1), there exist subgroups
$U'\le U\le G_{u_b}$ and $W'\le W\le G_{w_b}$  such that
we have
$$
\fs=(U\times W)/(U'\times W')=(U/U')\times (W/W'),
$$
It follows that $\fs=\fs_1\times\fs_2$, where $\fs_1=U/U'$ and $\fs_2=W/W'$. Assume that $\fs\in\fS_0(\cA)$. Let us verify that $\fs_1\preceq \fs_u$ for some $u\in V_1$; the fact that $\fs_2\preceq \fs_w$ for some $w\in V_2$, is proved similarly.

By \cite[Theorem~8.3]{Evdokimov2016}, we have $\cA_\fs=\mZ\,\fs$  and hence $\cA_{\fs_1}=\mZ\,\fs_1$. Without loss of generality, we may assume that $U\ne 1$ (otherwise one can take $u=1$). Then $G_{u_{i-1}}<U\le G_{u_i}$ for some $1\le i\le b$. Now if $U'<G_{u_{i-1}}$, then setting $\ft=G_{u_{i-1}}/U'$ and $\ft'=U/G_{u_{i-1}}$, we obtain
$$
\mZ\,\fs_1=\cA_{\fs_1}=\cA_{\ft}\wr\cA_{\ft'}\ne \mZ\,\fs_1,
$$
a contradiction. Thus, $G_{u_{i-1}}\le U'\le U\le G_{u_i}$, and we are done with $v=u_i$.
\eprf

In the following lemma, we show that  the inner $\cA$-multipliers are in a one-to-one correspondence with the elements of the direct product of the groups~$\aut(\fs_v)$, $v\in V$.

\lmml{030824s}
Put $\fS=\{\fs_v:\ v\in V\}$ and $\fM=\mult(\cA)$. Then
\nmrt
\tm{1} $\bar\fm\in\fM_\fS$ if and only if $\bar \fm=\{\sigma_v\in\aut(\fs_v):\ v\in V\}$,
\tm{2} the mapping  $\fM\to\fM_\fS$, $\fm\mapsto\fm_\fS,$ is a group monomorphism.
\enmrt
\elmm
\prf
The necessity in statement~(1) is trivial. To prove the sufficiency, we need to find an inner $\cA$-multiplier $\fm$ such that $\fm_\fS=\bar m$. First, put
\qtnl{040824a}
\fm_{\fs_{u,w}}=(\sigma_u,\sigma_w)
\eqtn
for all $u\in V_1$ and $w\in V_2$. Now let section $\fs\in\fS_0(\cA)$ be arbitrary. By Lemma~\ref{101123w}(2), there exist $u$ and $w$ such that $\fs\preceq \fs_{u,w}$. Put
\qtnl{040824b}
\fm_\fs:=(\fm_{\fs_{u,w}})^\fs.
\eqtn

We claim that $\fm_\fs$ does not depend on the choice of $u$ and $w$ for which  $\fs\preceq \fs_{u,w}$. Indeed, let $\fs\preceq \fs_{u',w'}$ for some $u'$ and $w'$ such that $u\ne u'$ or $w\ne w'$. Then either $u=u'$ and $w\ne w'$, or $u\ne u'$ and $w=w'$ (because the sections $\fs_{u,w}$ and $ \fs_{u',w'}$ have no common subsections if $u\ne u'$ and $w\ne w'$). Assume the first case occurs (in the second case the proof is similar). Then $\fs=\ft\times\ft'$ where $\ft\preceq\fs_u$ and $\ft'$ is a trivial subsection of both $\fs_w$ and~$\fs_{w'}$.  It follows that
$$
(\fm_{\fs_{u,w}})^\fs=({\sigma_u}^\ft, {\sigma_w}^{\ft'})=({\sigma_u}^\ft, \id_{\ft'})=({\sigma_{u'}}^\ft, {\sigma_{w'}}^{\ft'})=(\fm_{\fs_{u',w'}})^\fs,
$$
which proves the claim.

From the above, it follows that the  family $\fm$ defined by formulas~\eqref{040824a} and~\eqref{040824b}  satisfies the first condition in~\eqref{230824f} and the equality $\bar\fm=\fm_\fS$. To prove  that $\fm$ is an inner $\cA$-multiplier, we need to verify the second condition in~\eqref{230824f}. Let $\fs$ and~$\ft$ be projectively equivalent $\cA$-sections belonging to $\fS_0(\cA)$.  Then without loss of generality we may assume that both $\fs$ and~$\ft$ are subsections of a section $\fs_{u,w}$ for some~$u$ and~$w$.  The automorphism $\fm_{\fs_{u,w}}$ of the cyclic group $\fs_{u,w}$ acts by raising to a certain  integer power $\lambda$. Hence the automorphisms $\fm_\fs=(\fm_{\fs_{u,w}})^\fs$
and $\fm_\ft=(\fm_{\fs_{u,w}})^\ft$ act by raising to the powers $\lambda\pmod{|\fs|}$ and $\lambda\pmod{|\ft|}$, respectively. Since $|\fs|=|\ft|$ and the isomorphism $f_{\fs,\ft}:\fs\to\ft$ preserves raising to a power, we conclude that $ (\fm_\fs)^{f_{\fs,\ft}}=\fm_\ft$, as required.

To prove statement (2), assume that $\fm\in\fM$ is such that $\fm_\fs=\id_\fs$ for all $\fs\in\fS$. Note that if  $u\in V_1$ and $w\in V_2$, then $\fm_{\fs_{u,w}}=(\fm_{\fs_u},\fm_{\fs_w})=\id_{\fs_{u,w}}$.  Now let $\fs\in\fS_0(\cA)$. By Lemma~\ref{101123w}(2), we have $\fs\preceq \fs_{u,w}$ for some $u$ and $w$. By the first condition in~\eqref{230824f},  this implies that
$$
\fm_\fs=(\fm_{\fs_{u,w}})^\fs=(\id_{\fs_{u,w}})^\fs=\id_\fs.
$$
Thus the $\fm$ is the trivial inner $\cA$-multiplier and we are done.
\eprf

The basic sets of the S-ring $\cA$ can conveniently be divided into parts parameterized by the $\cA$-sets defined below; this parameterization will be used in Section~\ref{180524g}.   Namely, by Lemma~\ref{101123w}, for any two indices $1\le i,j\le b$, there is an $\cA$-subgroup $G_{i,j}=G_{u_i}\times G_{w_j}$. It follows that the difference
\qtnl{270724l}
D_{i,j}=G_{i,j}\setminus G_{i-1,j-1}
\eqtn
is an $\cA$-set. Since any nonidentity element of $G$ belongs to some of these differences, the union of all of them is equal to $G^\#$ and each nonidentity set $X\in\cS(\cA)$ is contained in at least one of them. Moreover, it is not hard to verify that
\qtnl{270724l1}
X\subseteq D_{i,j} \quad\Leftrightarrow\quad X=\pi_{i,j}^{-1}(x) \text{ for some } 1\ne x\in \im(\pi_{i,j}),
\eqtn
where $\pi_{i,j}:G_{i,j}\to G_{i,j}/G_{i-1,j-1}$ is the natural epimorphism.

\subsection{Inner $\cA$-multipliers associated with $\cG$}
Let $v=v_i$ be a vertex of the graph~$\cG$, $1\le i\le a$. Then $\fs_v=C_{n_i}$ is a (cyclic) group of order $n_i=p_iq_i$, where $p_i$ and~$q_i$ are distinct odd primes,  $p_i<q_i$. It follows that the group
$$
\aut(\fs_v)=\aut(C_{p_i})\times\aut(C_{q_i})
$$
has a uniquely determined Klein subgroup $K_v=\{k_v:\ k\in K\}$. For definiteness, we assume that for a fixed $k\in K^\#$ exactly one of the  conditions
$$
k_{v_i}\in\aut(C_{p_i})\qoq k_{v_i}\in\aut(C_{q_i})\qoq k_{v_i}\not\in\aut(C_{p_i})\text{ and }k_{v_i}\not\in\aut(C_{q_i})
$$
is satisfied simultaneously for  $i=1,\ldots,a$; in the latter case, $k_{v_i}$ is equal to the product of two involutions, one in $\aut(C_{p_i})$ and another in $\aut(C_{q_i})$.  Clearly, the mapping
\qtnl{150225t}
\pi_v:\bar K_v\to K_v,\ \bar k_v\mapsto k_v,
\eqtn
is a group isomorphism for all $v\in V$. The lemma below immediately follows from Lemma~\ref{030824s}.

\lmml{030824j}
Let $\fM=\mult(\cA)$ and $f\in K_V$. Then there is a unique inner $\cA$-multiplier $\fm_f\in\fM$ such that $(\fm_f)_{\fs_v}=\pi_v(f_v)$ for all $v\in V$, Moreover, the mapping
$$
\mu:K_V\to\fM,\ f\mapsto \fm_f,
$$
is a group monomorphism and $\im(\mu)=\{\fm\in\fM:\ \fm_{\fs_v}\in  K_v\}$.
\elmm

\section{The construction}\label{101024d}

\subsection{The main theorem}
In this section, we construct an S-ring $\cA^\star$ and a set of  $\cA^\star$-multipliers. They will be used to prove Theorem~\ref{090524d} in Section~\ref{180524g}. The relevant properties of~$\cA^\star$ are given in the theorem below.

\thrml{090524b}
Let $\cG\in\fB_a$ and $k$ a positive integer.  Assume that each separator of the graph $\cG$ has  cardinality greater than $k$. Then there exist
\nmrt
\tm{1} a totally normal circulant S-ring $\cA^\star$ of degree $n$ defined as in Theorem~\ref{090524d},  such that the scheme $\cX((\cA^\star)_0)$ is binary,
\tm{2} an outer $\cA^\star$-multiplier $\fc$ such that $\iso(\cA^\star,\varphi_\fc)=\varnothing$, and 
\tm{3} a $\lfloor  k/6\rfloor$-local system of inner $\cA^\star$-multipliers with respect to  $\fc$.
\enmrt
\ethrm

The proof will be given in the end of the section. In what follows, we keep notation of Section~\ref{190524h}.

\subsection{Construction of the S-ring $\cA^\star$}\label{070624a}
Let $\cX$  be the coherent configuration defined in Subsection~\ref {121123u}. By the condition~(K3), we have $\aut(\cX)\le K_V$.  Let
$$
\fM=\mu(\aut(\cX)),
$$
where $\mu$ is the monomorphism defined in Lemma~\ref{030824j}. As a consequence of  conditions~(K1), (K3), and (K4), and Lemma~\ref{030824j}, we have the following statement, in which  we abbreviate $\fM_v:=\fM_{\fs_v}$ and $\fM_{v,v'}:=\fM_{\{\fs_{v^{}},\fs_{v'}\}}$ for all $v,v'\in V$, and denote by $\rho_v:K_{v^{}}\to \sym(2)$  the natural epimorphism with kernel $\grp{k_v}$.

\lmml{070724a}
In the above notation, $\fM_v=K^{}_v$ for all $v\in V$. Furthermore, if $v'\in V$, then
$$
\fM_{v,v'}=\css
\{(x,y)\in K_{v^{}}\times K_{v'}:\ \rho_{v^{}}(x)=\rho_{v'}(y)\} &\text{if $(v,v')\in E$,}\\
K_{v^{}}\times K_{v'}  &\text{otherwise.}\\
\ecss
$$
\elmm

Let $\cA$ be the S-ring defined in Subsection~\ref{130524f}. Denote by  $\cA^\star$   the algebraic fusion of  $\cA$ with respect  to the subgroup~$\fM\le \fM(\cA)$,
$$
\cA^\star=\cA^\fM.
$$
Since $\cA$ and $\cA^\star$ have the same $\cA$-subgroups (formula~\eqref{040824q}), the subsets $D_{i,j}$ defined by formula~\eqref{270724l} are $\cA^\star$-sets for all~$i,j$.

Let  $\fs\preceq\fs_{u,w}$ for some  $u\in V_1$ and $w\in V_2$. Then from Lemma~\ref{070724a}, it follows that
$$
\fM_\fs\le (K_u)^\fs\times (K_w)^\fs \le\aut(\fs).
$$
Since the group $\fs$ is of odd order whereas the groups $(K_u)^\fs$ and $(K_w)^\fs$ are of order at most~$4$, the cyclotomic  S-ring $\cyc(\fM_\fs,\fs)$ is normal by \cite[Theorem~6.1]{EvdP2003}. By Lemma~\ref{120724a}, this implies that
\qtnl{130724a}
\fS_0(\cA^\star)=\fS^{}_0(\cA)\qaq (\cA^\star)_0=\cA.
\eqtn

\lmml{120724d}
 The S-ring $\cA^\star$ is totally normal.
\elmm
\prf
Let  $\fs\in\fS_0(\cA^\star)$. Then $\fs\in \fS^{}_0(\cA)$ by the left-hand side of formula~\eqref{130724a}. By Lemma~\ref{101123w}(2), this implies that $\fs\preceq\fs_{u,w}$ for some $u\in V_1$ and $w\in V_2$. Since by the above,  the S-ring $(\cA^\star)_\fs=\cyc(\fM_\fs,\fs)$ is normal, we are done.
\eprf

\subsection{An outer multiplier}\label{140524e}
In the notation of Subsection~\ref{121123u}, fix an arbitrary edge $e_0=\{u_{i_0},w_{j_0}\}$ of the graph $\cG$ and the algebraic automorphism  $\psi=\psi_{e_0}$ of the coherent configuration~$\cX$.  From condition~(K4), it follows that the group $\aut(\cX)^{\Omega_{e_0}}$ has index~$2$ in the group $\bar K_{u_{i_0}}\times\bar K_{w_{j_0}}$. Choose any permutation $\bar k_0$ of this group, not belonging to $\aut(\cX)^{\Omega_{e_0}}$, and define the permutation $k_0\in\aut(\fs_{e_0})$ by  the conditions
$$
(k_0)_{\fs_{u_{i_0}}}=\pi_{u_{i_0}}((\bar k_0)_{u_{i_0}})\qaq
(k_0)_{\fs_{w_{j_0}}}=\pi_{w_{j_0}}((\bar k_0)_{w_{j_0}}).
$$
The definition of the monomorphism $\mu$ implies that if $f\in K_V$ and $f^{e_0}=\bar k_0$, then $(\mu(f))_{e_0}=k_0$.

Let $u\in V_1$ and $w\in V_2$. Then the section $\fs=\fs_{u,w}$ belongs to $\fS_0(\cA)$, see Lemma~\ref{101123w}(2). By formula~\eqref{130724a}, this implies that $\fs\in\fS_0(\cA^\star)$. Let us define  a coset $\fc_\fs$ of the subgroup~$\fM_\fs\le\aut(\fs)$ by the formula
\qtnl{100624a}
\fc_\fs=\css
k_0\cdot\fM_\fs  & \text{if $\fs=\fs_{e_0}$,}\\
\fM_\fs & \text{otherwise,}\\
\ecss
\eqtn
where $\fs_{e_0}=\fs_{u_{i_0},w_{j_0}}$. We extend this definition to all sections $\fs\in\fS_0(\cA^\star)$ by setting $\fc_\fs=\fc_{\fs_{u,w}}$ for suitable $u,w$ (see Lemma~\ref{101123w}(2) and formula~\eqref{130724a}). Furthermore, by formula~\eqref{040824g}, we have
$$
(\cA^\star)_\fs=\cyc(\fM_\fs,\fs).
$$
This S-ring is normal by Lemma~\ref{120724d} and hence $\fM_\fs=\aut_{\cA^\star}(\fs)$ by Lemma~\ref{280724b}(3). Using arguments in the proof of Lemma~\ref{030824s}, one can verify that the family $\fc:=\fc(e_0)$ of the cosets $\fc_\fs$, $\fs\in \fS_0(\cA^\star)$, is an outer $\cA^\star$-multiplier.

In the following statement, we establish a relationship between the isomorphisms of the coherent configuration $\cX$, the corresponding inner multipliers in $\fM$, and the outer multiplier~$\fc$.

\lmml{140524a}
Let   $f\in K_V$ and $e=\{v,v'\}$, where $v,v'\in V$ and $v\ne v'$. Then
\qtnl{100624w}
f^{\Omega_e}\in \iso(\cX_{\Omega_e},\psi_{\Omega_e})\quad\Leftrightarrow\quad \text{either }e\not\in E\text{ or }  (\fm_f)_e\in \fc_e,
\eqtn
where $\fm_f=\mu(f)$ with $\mu$ as in Lemma~{\rm \ref{030824j}},  $(\fm_f)_e=(\fm_f)_{\fs_e}$, and $\fc_e=\fc_{\fs_e}$.
\elmm
\prf
We note that if  $e\not\in E$, then $\iso(\cX_{\Omega_e},\psi_{\Omega_e})=\bar K_{v^{}}\times \bar K_{v'}$  and $f^{\Omega_e}\in \iso(\cX_{\Omega_e},\psi_{\Omega_e})$ by  condition~(K5). Assume that $e\in E$.

Let $e\ne e_0$. Then  $\iso(\cX_{\Omega_e},\psi_{\Omega_e})=\aut(\cX_{\Omega_e})$ by condition~(K4). Therefore, $f^{\Omega_e}\in \iso(\cX_{\Omega_e},\psi_{\Omega_e})$ if and only if
$$
(\fm_f)_e\in\mu_e(\aut(\cX_{\Omega_e}))=\fM_e=\fM_{\fs_{v,v'}}=\fc_{\fs_{v,v'}}=\fc_e,
$$
see Lemma~\ref{070724a}, where $\mu_e$ is the restriction of the monomorphism~$\mu$ to~$e$.

Finally, let $e=e_0$, i.e., $\{v,v'\}=\{u_{i_0},w_{j_0}\}$. Then by condition~(K4), the group $\aut(\cX_{\Omega_e})$ has index~$2$ in the group $\bar K_{u_{i_0}} \times\bar K_{w_{j_0}}$. By virtue of condition~(K5), we have
$$
\iso(\cX_{\Omega_e},\psi_{\Omega_e})=\bar k_0\cdot\aut(\cX_{\Omega_e}),
$$
where $\bar k_0$ is the permutation chosen at the beginning of this subsection. Thus,
$f^{\Omega_e}\in \iso(\cX_{\Omega_e},\psi_{\Omega_e})$ if and only if
$$
(\fm_f)_e\in\mu_e(\bar k_0\cdot  \aut(\cX_{\Omega_{e_0}}))=\fM_{e_0}=\fM_{\fs_{u_{i_0},w_{j_0}}}=\fc_{\fs_{u_{i_0},w_{j_0}}}=\fc_{e_0},
$$
as required
\eprf

\subsection{Proof of Theorem~\ref{090524b}} Let $\cA^\star$ be the  S-ring over the group~$G$, defined in Subsection~\ref{070624a}. By Lemma~\ref{120724d}, $\cA^\star$ is totally normal. By formula~\eqref{130724a}, we have ~$(\cA^\star)_0=\cA$. Hence the scheme $\cX((\cA^\star)_0)$ is binary by Lemma~\ref{250824b}. Thus it suffices to verify that statements~(2) and~(3) of the theorem in question hold true for the outer $\cA^\star$-multiplier $\fc$ defined in Subsection~\ref{140524e}.

To prove statement~(2), assume on the contrary that there is  $g\in\iso(\cA^\star,\varphi_\fc)$. Without loss of generality, we may assume that $g$ is a normalized isomorphism.
Since $G$ is of odd order and $\cA^\star$ is totally normal, the family
$$
\fm=\{g^\fs:\ \fs\in\fS_0(\cA^\star)\}
$$
is an inner $\cA^\star$-multiplier, see the end of Subsection~\ref{180225a}. Furthermore, $\fm$ is an inner $\cA$-multiplier by the first equality in~\eqref{130724a}. In addition, $\fm\le\fc$ and hence
$$
\fm_{\fs_v}=g^{\fs_v}\in\fc_{\fs_v}=\fM_{\fs_v}=K_v\qquad v\in V.
$$
By Lemma~\ref{030824j}, this shows  that $f=\mu^{-1}(\fm)$ is an  isomorphism of the coherent configuration~$\cX$, belonging to~$K_V$. Besides, let $e=\{v,v'\}$, where  $v,v'\in V$ and $v\ne v'$. Then
$$
(\fm_f)_e=(\pi_{v^{}}(f_{v^{}}),\pi_{v'}(f_{v'}))=(g_{v^{}},g_{v'})=\fm_e\in \fc_e,
$$
where $\pi_{v^{}}$ and $\pi_{v'}$ are the isomorphisms defined in formula~\eqref{150225t}.  By  Lemma~\ref{140524a}, this implies  that $f^{\Omega_e}\in\iso(\cX_{\Omega_e},\psi_{\Omega_e})$. It follows that $f\in\iso(\cX,\psi)$, which contradicts Lemma~\ref{050524t}.

Let us prove statement~(3).  Put $m=\lfloor  k/6\rfloor$. The graph $\cG$ has no separators of cardinality at most~$3\cdot (2m)\le k$.  By Lemma~\ref{050524t},  there exists an $(2m)$-local system
\qtnl{070624r}
\{f_\Delta\in K_V: \Delta \in F_{2m}(\cX)\}
\eqtn
of isomorphisms for the coherent configuration $\cX$ with respect to  the algebraic isomorphism $\psi$.

Let $\fS\subseteq \fS_m$, where  $\fS_m=\fS_m(\cA^\star)$. By formula in~\eqref{130724a} and Lemma~\ref{101123w}(2), there is a minimal (by inclusion) subset $\bar V=\bar V(\fS)$ of the set~$V$ such that if $\fs\in\fS$, then $\fs\preceq\fs_{u,w}$ for some $u\in V_1\cap\bar V$ and $w\in V_2\cap \bar  V$. Clearly,
$$
|\bar  V|\le 2\cdot |\fS|\le 2\cdot m
$$
and hence $\Omega_{\bar V}\in F_{2m}(\cX)$.  Thus the $(2m)$-local system~\eqref{070624r} contains an isomorphism $f_\Delta\in K_V$ of the coherent configuration~$\cX$, where $\Delta=\Omega_{\bar V}$. Denote by  $\fm(\fS)$ the inner $\cA$-multiplier $\mu(f_\Delta)$.

To complete the proof, it suffices to verify that the family
$$
\{\fm(\fS)\in\fM(\cA^\star):\ \fS\in \fS_m\}
$$
is an $m$-local system of inner $\cA^\star$-multipliers with respect to  the outer $\cA^\star$-multiplier~$\fc$.

To verify the first relation in~\eqref{080423a}, let $\fS\in \fS_m$ and $\fs\in \fS$. 
By Lemma~\ref{101123w}(2), we may assume without loss of generality that $\fs=\fs_{u,w}$ for some vertices $u\in V_1\cap\bar V$ and $w\in V_2\cap \bar  V$, where $\bar V=\bar V(\fS)$. Since $u,w\in \bar V$, we have  $\Omega_{u,w}\subseteq\Omega_{\bar V}:=\Delta$.  By the definition of~$f_\Delta$, this implies that $(f_\Delta)^{\Omega_e}\in \iso(\cX_{\Omega_e},\psi_{\Omega_e})$, where $e=\{u,w\}$. By Lemma~\ref{140524a}, either  $e\not\in E$, or $e\in E$ and
$$
\fm(\fS)_\fs=\mu(f_\Delta)_\fs=\mu(f_\Delta)_{\fs_e}\in\fc_e=\fc_\fs.
$$
Finally, if $e\not\in E$, then  $\fm(\fS)_\fs\in K_u\times K_w=\fc_\fs$.

To verify the second relation in~\eqref{080423a},  let $\fS,\fS'\in\fS_m$.  Put $\Delta=\Omega_{\bar V^{}}$ and $\Delta'=\Omega_{\bar V'}$, where $\bar V=\bar V(\fS)$ and $\bar V'=\bar V(\fS')$. As above, we have $\Delta,\Delta'\in F_{2m}(\cX)$. By the definition of  $(2m)$-local system, this implies that
$$
(f_\Delta {f_{\Delta'}}^{-1})^{\Delta\cap\Delta'}=(h_{\Delta,\Delta'})^{\Delta\cap\Delta'}
$$
for some $h_{\Delta,\Delta'}\in\aut(\cX)$. It follows that $f_\Delta {f_{\Delta'}}^{-1}=f_{\Delta,\Delta'}h_{\Delta,\Delta'}$, where the isomorphism $f_{\Delta,\Delta'}:=f_\Delta {f_{\Delta'}}^{-1}h_{\Delta,\Delta'}^{-1}$ belongs to $K_V$ and is identical on $\Delta\cap\Delta'$. Applying the monomorphism~$\mu$ from Lemma~\ref{030824j} to the last equality, we obtain
\qtnl{090624a}
\fm(\fS)\fm(\fS')^{-1}=\mu(f_{\Delta,\Delta'}h_{\Delta,\Delta'})=\mu(f_{\Delta,\Delta'})\mu(h_{\Delta,\Delta'}).
\eqtn
By Lemma~\ref{120724a}, the inner $\cA^\star$-multiplier  $\fm(\fS,\fS'):=\mu(h_{\Delta,\Delta'})$ belongs to the group  $\fM\le \fM_{aut}(\cA^\star)$. Moreover, the inner $\cA^\star$-multiplier $\mu(f_{\Delta,\Delta'})$ is identical on $\fS\cap \fS'$. Thus, taking restriction of the equality \eqref{090624a} to $\fS\cap \fS'$, we get
$$
\fm(\fS)_{\fS\cap \fS'}=\fm(\fS')_{\fS\cap \fS'}(\fm(\fS,\fS'))_{\fS\cap \fS'},
$$
as required.

\section{Proof of Theorem~\ref{090524d} and Corollary~\ref{250924a}}\label{180524g}
\subsection{Reduction from S-rings to graphs} Let $\cG\in\fB_a$ be an  $\varepsilon$-expander for some real $\varepsilon>0$. This means that for every vertex subset $\Delta$ containing at most half of the vertices of $\cG$, we have $|\partial \Delta|\ge \varepsilon |\Delta|$, where $\partial \Delta$ is the set of all vertices outside $\Delta$, each adjacent to at least one vertex of~$\Delta$. It is not hard to prove (see, e.g., \cite[Lemma 4.2.7]{CP2019}) that each separator of the graph~$\cG$  has  cardinality greater than $k=\lfloor c_\varepsilon\cdot 2a\rfloor+1$, where $c_\varepsilon ={\varepsilon\over{4+\varepsilon}}$.

Let $\cA^\star$ be the totally normal S-ring over a cyclic  group~$G$ of order~$n$ and $\varphi^\star=\varphi_\fc$ the algebraic isomorphism of $\cA^\star$, that were defined in Theorem~\ref{090524b}; in particular,  $\cA^\star$ admits  a $\lfloor k/6\rfloor$-local system of inner $\cA^\star$-multipliers with respect to $\fc$ and $\iso(\cA^\star,\varphi^\star)=\varnothing$. Then
\qtnl{270424r}
\iso(\cX^\star,\varphi^\star)=\varnothing,
\eqtn
where $\cX^\star=\cX(\cA^\star)$. Furthermore, from Theorem~\ref{310324ab}, it follows that $\cX^\star$  is $\WL_m$-equivalent to itself with respect to the algebraic isomorphism $\varphi^\star$, where $m=\lfloor \sqrt{\lfloor k/6\rfloor}-1\rfloor=c\cdot\sqrt{a}$ for a certain real constant $0<c<1$. This fact together with relation~\eqref{270424r} shows that
\qtnl{270724n}
\dim_{\scriptscriptstyle\WL}(\cX^\star)\ge c\cdot\sqrt{a},
\eqtn
see Section~\ref{040924e}. In the rest of this section, we construct a circulant graph $\cG^\star$ over the group $G$ and verify  that
\qtnl{300824a}
\WL(\cG^\star)=\cX^\star.
\eqtn
This is enough to complete the proof. Indeed, in accordance with formulas~\eqref{040924a} and~\eqref{270724n}, we then have $\dim_{\scriptscriptstyle\WL}(\cG^\star)=\dim_{\scriptscriptstyle\WL}(\cX^\star)\ge c\cdot\sqrt{a}$, as required.

\subsection{The definition of the graph $\cG^\star$}\label{091024d} To construct the graph~$\cG^\star$, for any two indices $1\le i,j\le b$, we choose an arbitrary basic set $X_{i,j}\in \cS(\cA^\star)$ contained in the difference~$D_{i,j}$ defined by formula~\eqref{270724l} and such that
 $(X_{i,j})_{\fs_{i,j}}$ is a highest basic set of the S-ring $\cA_{\fs_{i,j}}$, where $\fs_{i,j}:=\fs_{u_i,w_j}$.  We have
\qtnl{290724c}
|X_{i,j}|=|\pi_{i,j}^{-1}((X_{i,j})_{\fs_{i,j}})|=|G_{i-1,j-1}|\cdot\css
4 & \text{if $i=1$ or $j=1$},\\
8 & \text{if $i,j>1$ and $\{u_i,w_j\}\in E$},\\
16 & \text{otherwise},\\
\ecss
\eqtn
where $G_{i-1,j-1}$ and $\pi_{i,j}$ are as  in formulas~\eqref{270724l} and~\eqref{270724l1}, respectively. Denote by~$X^\star$ the union of all these basic sets~$X_{i,j}$. Then the Cayley graph
$$
\cG^\star=\cay(G,X^\star)
$$
is  circulant of order~$n$.

\subsection{Proof of equality~\eqref{300824a}}
The scheme $\WL(\cG^\star)$ is  a Cayley scheme over the group $G$. It is associated with the  S-ring $\cA'=\WL_G(\und{X^\star})$. Since $X^\star$ is also an $\cA^\star$-set, we have $\cA^\star\ge \cA'$. It remains to prove the reverse inclusion
\qtnl{300824r}
\cA^\star\le \cA',
\eqtn
i.e., that each basic set of $\cA^\star$ is an $\cA'$-set.  To do this, we need several auxiliary statements.

\clml{280724g1}
$G_{u_b}$ and $G_{w_b}$ are $\cA'$-groups.
\eclm
\prf
We prove the statement for $G_{u_b}$ only; the statement for $G_{w_b}$ is proved similarly.
Recall that  the number $n_1$ is a product of two distinct primes, say $p_1$ and $q_1$. We apply Theorem~\ref{261009b}(2) for $\cA=\cA'$ and $X=X^\star$. Note that if $x\in X^\star$ and $H$ the subgroup defined before that theorem, then $xH\subseteq X^\star$ unless  $x$  belongs to the union of all~$X_{i,j}$ with $i=0,1$.  Therefore,
$$
((X^\star)^{[p_1]})^{[q_1]}=\Bigl(\bigcup_{i,j=0}^b(X_{i,j})^{[p_1]}\Bigl)^{[q_1]}=
\Bigl(\bigcup_{j=0}^b(X'_{0,j}\cup X'_{1,j})\Bigr)^{[q_1]}=\bigcup_{j=0}^b X''_{0,j}
$$
where for all $i,j$, we set $X'_{i,j}=\sum_{x\in X_{i,j}} x^{p_1}$ and $X''_{i,j}=\sum_{x\in X'_{i,j}} x^{q_1}$. The  right-hand side is contained in $G_{u_b}$ and contains the set $X''_{0,b}$ generating $G_{u_b}$. Thus, $G_{u_b}$  is an $\cA'$-group.
\eprf

\clml{280724g}
$G_{u_i}$ and $G_{w_i}$ are $\cA'$-groups, $i=1,\ldots,b$.
\eclm
\prf
From Claim~\ref{280724g1}, it follows that $X:=X^\star\cap G_{u_b}$ is an $\cA'$-set. Note that $X=X_{0,0}\cup X_{1,0}\cup\ldots\cup X_{b,0}$ and $X_{i,0}\subseteq D_{i,0}$ for  $i=0,1,\ldots,b$. Since each $X_{i,0}$ is a basic set of the S-ring $\cA^\star$, we have
\qtnl{290724a}
(\und{X}\cdot\und{X})\circ\und{X}=\sum_{i=0}^ba_i\und{X_{i,0}}
\eqtn
for some nonnegative integers  $a_i$, where $\circ$ denotes the componentwise product in~$\mZ G$.  The integer $a_i$ is equal to the number of pairs $(y,z)\in X\times X$ such that $yz=x$ for an arbitrary  fixed $x\in X_{i,0}$.

Let $y\in D_{i',0}$,  where $0\le i'\le b$. Assume that  $i'>i$. Then  $y^{-1}x\in y^{-1}G_{i,0}\subseteq X$, and $x=yz$ with $z:=y^{-1} x$.  Similarly, if $i'<i$, then $x=yz$ for some $z\in X$. Now suppose that $i'=i$. In this case, if $i>0$, then the set $(X_{\fs_{i,0}})\cdot  (X_{\fs_{i,0}})$ contains no element of $X_{\fs_{i,0}}$ (this follows from the assumption that $|\fs_{i,0}|\ge 5$). Thus,
$$
a_i=|X|-|X_{i,0}|, \quad i=1,\ldots b.
$$	
Besides, $a_0=|X|$.  Together with  formula~\eqref{290724c}, this shows that the integers  $a_0,a_1.\ldots,a_b$ are pairwise distinct. It easily follows that $X_{i,0}$ is an $\cA'$-set for all~$i$. Therefore, $G_{u_i}=\grp{X_{i,0}}$ is an $\cA'$-group, $i=1,\ldots,b$. In a similar way, one can prove that so is the group $G_{w_j}=\grp{X_{0,j}}$, $j=1,\ldots,b$.
\eprf

\crllrl{aa}
Let $0\le i,j\le b-1$. Then $\fs_{i,j}\in\fS(\cA')$ and $D_{i,j},X_{i,j}\in\cS(\cA')^\cup$.
\ecrllr
\prf
By Claim~\ref{280724g}, $G_{i,j}=G_{u_i}\times G_{w_j}$ is an $\cA'$-group for all $i,j$. Then obviously $\fs_{i,j}$ is an $\cA'$-section and $D_{i,j}$ is an $\cA'$-set. Moreover, so is the subset
$X_{i,j}=X^\star\cap D_{i,j}$.
\eprf

Let $X$ be a nonidentity basic set of $\cA^\star$. Then there are indices $0\le i,j\le b-1$ such that $X\subseteq D_{i,j}$. Put $\fs=\fs_{i,j}$. From Lemma~\ref{120724d}, it follows that  $(\cA^\star)_\fs$ is a  normal S-ring over $\fs$. By the choice of the basic set $X_{i,j}$, Lemma~\ref{280724b} implies that $(\cA^\star)_\fs$ is the smallest S-ring over~$\fs$ that contains $(X_{i,j})_\fs$. By Corollary~\ref {aa}, this implies that $(\cA')_\fs\ge (\cA^\star)_\fs$. Therefore,
$$
X^{}_\fs\in\cS(\cA^\star_\fs)^\cup\subseteq  \cS(\cA'_\fs)^\cup.
$$
By  formula~\eqref{270724l1}, this shows that  $X=\pi_{i,j}^{-1}(X_\fs)$ is an $\cA'$-set, which proves inclusion~\eqref{300824r}.

\subsection{Proof of Corollary~\ref{250924a}}
Recall (see, e.g., \cite{LubotzkyPS1986} and references therein) that a $k$-regular graph $\cG$ is \emph{Ramanujan} if all of its nontrivial eigenvalues lie between $-2\sqrt{k-1}$ and $2\sqrt{k-1}$.  Any such graph is an $\varepsilon$-expander with $2\varepsilon=1-\lambda/k$, where~$\lambda$ is the second largest (in absolute value) eigenvalue of the adjacency matrix of~$\cG$. By Theorem~\ref{090524d}, it suffices to find an infinite family of bipartite cubic $3$-connected Ramanujan graphs.

Let $d\in\mN$ be even. Denote by $\cG$ a connected cubic Ramanujan graph on  $2^{3d}-2^d$ vertices, constructed in~\cite[Theorem~5.13]{Morgenstern1994}. Note that it  is a Cayley graph of $\PSL(2,2^d)$.  Let  $\cG_d$ be the $2$-cover of $\cG$, i.e., the vertex set of $\cG_d$ is the disjoint union of the vertex sets of two copies of~$\cG$ and two vertices  of $\cG_d$ are adjacent if they lie in different copies and the corresponding vertices of $\cG$ are adjacent. Obviously, the graph $\cG_d$ is cubic and bipartite. Moreover, it is a Ramanujan graph (see, e.g., a remark in~\cite[p.~311]{MarcusSS2015}).

It remains to verify that the graph $\cG_d$ is $3$-connected. From the fact that $\cG$ is a Cayley graph, it easily follows that the graph $\cG_d$ is vertex transitive. The connectivity of it can be equal neither $1$ (because it contains at least one vertex which is not a cut point), nor $2$ (because any vertex transitive graph of connectivity~$2$ is a polygon~\cite[Corollary~3A]{Wat1970}). Thus  the (cubic) graph $\cG_d$ is $3$-connected.

\section*{Acknowledgment}
The authors would like to thank Gang Chen for his useful discussion and helpful suggestions.

\bibliographystyle{amsplain}

\end{document}